\documentclass[]{interact}
\usepackage{dsfont}
\usepackage[caption=false]{subfig}

\usepackage[numbers,sort&compress]{natbib}
\bibpunct[, ]{[}{]}{,}{n}{,}{,}
\renewcommand\bibfont{\fontsize{10}{12}\selectfont}
 
\theoremstyle{plain}
\newtheorem{theorem}{Theorem}[section]
\newtheorem{lemma}[theorem]{Lemma}

\newtheorem{assumption}[theorem]{Assumption}

\theoremstyle{definition}
\newtheorem{definition}[theorem]{Definition}

\theoremstyle{remark}
\newtheorem{remark}{Remark}

\usepackage{graphicx}
\usepackage[dvipsnames]{xcolor}

\usepackage{algorithm}
\usepackage{algorithmic}
\usepackage{caption}

\usepackage{mathtools}
\usepackage{amsfonts}
\usepackage{amssymb}
\usepackage{url}
\usepackage{dirtytalk}
\usepackage{tikz}
\usepackage{nicefrac}
\usepackage{hyperref}
\usepackage{multirow}
\usepackage{cancel}

\def\rev#1{{\color{black}#1}}

\begin{document}

\articletype{ARTICLE TEMPLATE}

\title{Gradient-Type Methods For Decentralized Optimization Problems With Polyak-{\L}ojasiewicz Condition Over Time-Varying Networks  \footnote{This work was supported by a grant for research centers in the field of artificial intelligence, provided by the Analytical Center for the Government of the Russian Federation in accordance with the subsidy agreement (agreement identifier 000000D730321P5Q0002) and the agreement with the Moscow Institute of Physics and Technology dated November 1, 2021 No. 70-2021-00138.
}
}

\author{
\name{Ilya Kuruzov\textsuperscript{a}\thanks{CONTACT Ilya Kuruzov. Email: kuruzov.ia@phystech.edu } and Mohammad Alkousa\textsuperscript{a,b} and Fedor Stonyakin\textsuperscript{a,c} and Alexander Gasnikov\textsuperscript{a,b,d}}
\affil{
\textsuperscript{a}Moscow Institute of Physics and Technology, Dolgoprudny, Russia;
\textsuperscript{b}National Research University Higher School of Economics, Moscow, Russia;
\textsuperscript{c}V. Vernadsky Crimean Federal University, Simferopol, Republic of Crimea; \textsuperscript{d}ISP RAS Research Center for Trusted Artificial Intelligence
}
}
\maketitle

\begin{abstract}
This paper focuses on the decentralized optimization (minimization and saddle point) problems with objective functions that satisfy Polyak-{\L}ojasiewicz condition (PL-condition). The first part of the paper is devoted to the minimization problem of the sum-type cost functions. In order to solve a such class of problems, we propose a gradient descent type method with a consensus projection procedure and the inexact gradient of the objectives. Next, in the second part, we study the saddle-point problem (SPP) with a structure of the sum, with objectives satisfying the two-sided PL-condition. To solve such SPP, we propose a generalization of the Multi-step Gradient Descent Ascent method with a consensus procedure, and inexact gradients of the objective function with respect to both variables. Finally, we present some of the numerical experiments, to show the efficiency of the proposed algorithm for the robust least squares problem. 
\end{abstract}


\section{Introduction}\label{sec:intro}

In this paper, firstly we study a sum-type minimization problem
\begin{equation}\label{main_sum_problem}
\min_{x \in \mathbb{R}^d} \left\{f(x) : = \frac{1}{n}\sum_{i=1}^{n} f_i(x)\right\},
\end{equation}
where the functions $f_i$ are generally non-convex and stored separately by nodes in a communication network, which is represented by a non-directed graph $\mathcal{G}= (E, V)$. The graph $\mathcal{G}$, possibly, can have a change over time structure. This problem, where the object function depends on distributed data is typed as a decentralized optimization problem. We assume that $f$ satisfies the well-known Polyak-{\L}ojasiewicz condition (for brevity, we write PL-condition). This condition was originally introduced by Polyak \cite{polyak_1963}, who proved that it is sufficient to show the global linear convergence rate for the gradient descent without assuming convexity. The PL-condition is very well studied by many researchers in many different works for many different settings of optimization problems and has been theoretically verified for objective functions of optimization problems arising in many practical problems. For example, it has been proven to be true for objectives of over-parameterized deep networks \cite{du2019gradient}, learning LQR models \cite{fazel2018global}, phase retrieval \cite{sun2018geometric}, Generative adversarial imitation learning of linear quadratic (see Example 3.1 in \cite{Baraz_non_convex_minimax}). More discussions of PL-condition and many other simple problems can be found in \cite{karimi2016linear}. 

This type of problems arises in different areas: distributed machine learning \cite{Mlbase}, resource allocation problem \cite{resource_allocation_problem} and power system control \cite{Nedic_power_control}. 


The problem \eqref{main_sum_problem} can be reformulated as a problem with linear constraints. For this, let us assign each agent in the network a personal copy of parameter vector $x_i \in \mathbb{R}^d$ (column vector) and introduce
\begin{equation}\label{function_FX}
\mathbf{X}:=\left(x_{1}^{\top} \, \ldots \; x_{n}^{\top}\right)^{\top} \in \mathbb{R}^{n \times d}, \quad  F(\mathbf{X})=\sum_{i=1}^{n} f_{i}\left(x_{i}\right).
\end{equation}


Now we equivalently rewrite problem \eqref{main_sum_problem} as
\begin{equation}\label{equiv_main_sum_problem}
\min_{\mathbf{X} \in \mathbb{R}^{n \times d}} F(\mathbf{X}) = \sum_{i=1}^{n}f_i(x_i),\quad \text{s.t.} \;\; x_1 = \dots= x_n,
\end{equation}
where it has the same optimal value as problem \eqref{main_sum_problem}. This reformulation increases the number of variables but induces additional constraints at the same time. 

Let us denote the set of consensus constraints $\mathcal{C} = \{\mathbf{X}| x_1 = \dots = x_n\}$. Also, for each $\mathbf{X} \in \mathbb{R}^{n \times d}$ denote the average of its columns $\overline{x} = \frac{1}{n} \sum_{i=1}^{n}x_i \in \mathbb{R}^d$ and introduce its projection onto constraint set
\[
\overline{\mathbf{X}}=\frac{1}{n} \mathbf{1}_{n} \mathbf{1}_{n}^{\top} \mathbf{X}= \Pi_{\mathcal{C}}(\mathbf{X})=\left(\overline{x}^{\top}\, \ldots \; \overline{x}^{\top}\right)^{\top}  \in \mathbb{R}^{n \times d}. 
\]

Note that $\mathcal{C}$ is a linear subspace in $\mathbb{R}^{n \times d}$, and therefore projection operator $\Pi_{\mathcal{C}}$ is linear. 

Secondly, we study the saddle-point problem with a structure of the sum
\begin{equation}\label{base_saddle_point_problem}
\min_{x\in \mathbb{R}^{d_x}}\max_{y\in \mathbb{R}^{d_y}}\left\{ \phi(x,y)=\frac{1}{n}\sum\limits_{i=1}^n \phi_i(x,y)\right\},
\end{equation}
where functions $\phi_i(\cdot, \cdot)$  are non-convex (with respect to the variable $x$ for each $y$) and smooth (i.e., with Lipschitz-continuous gradient). These functions \rev{can be calculated only} separately in different nodes in the communication network, which is represented by a non-directed graph $\mathcal{G}$. 

Problems of type \eqref{base_saddle_point_problem} arise in many applications, such that Generative adversarial network \cite{Goodfellow}, adversarial training \cite{Goodfellow_scale}, and fair training \cite{Baraz_non_convex_minimax}. In our work with problem \eqref{base_saddle_point_problem}, we assume that  $\phi(\cdot, y)$ and $-\phi(x,\cdot)$ satisfy the PL-condition (see Assumption \eqref{PL_assump}, below).

Let $y(x) : = \arg\max_{y\in \mathbb{R}^{d_y}} \phi(x,y)$, and let we set $f_i(x) : = \phi_i(x,y(x))$, then problem \eqref{base_saddle_point_problem} can be rewritten in the form of the minimization problem \eqref{main_sum_problem}, i.e.,
\begin{equation}\label{base_x_problem}
\min_{x\in \mathbb{R}^{d_x}}\left\{\max_{y\in \mathbb{R}^{d_y}} \left\{\phi(x,y)=\frac{1}{n}\sum\limits_{i=1}^n \phi_i(x,y)\right\}\right\} = \min_{x \in \mathbb{R}^{d_x}} \left\{f(x) : = \frac{1}{n}\sum_{i=1}^{n} f_i(x)\right\}.
\end{equation}

Let  $\mathcal{C}_x := \left\{\mathbf{X} \in \mathbb{R}^{n \times d_x}|x_1=\ldots = x_n\right\}$,  $\mathcal{C}_y: =\left \{\mathbf{Y} \in \mathbb{R}^{n \times d_y}|y_1=\ldots = y_n\right\}$ and  $\Phi(\mathbf{X},\mathbf{Y}) := \sum\limits_{i=1}^n \phi_i(x_i,y_i)$. So, we can rewrite \eqref{base_saddle_point_problem}, in the similar way like \cite{rogozin2021accelerated,rogozin_accel_do},  in the following form
\begin{equation}\label{base_saddle_point_problem_G}
\min_{\mathbf{X}\in \mathcal{C}_x}\max_{\mathbf{Y}\in \mathcal{C}_y} \Phi(\mathbf{X},\mathbf{Y}).
\end{equation}

Similarly to the first proposed minimization problem \eqref{main_sum_problem}, note that the problem \eqref{base_saddle_point_problem}, also can be rewritten in the same form
\begin{equation}\label{base_big_x_problem}
\min_{\mathbf{X}\in \mathcal{C}_x} F(\mathbf{X}) = \sum_{i=1}^{n}f_i(x_i).
\end{equation}

So, for solving the problem \eqref{base_x_problem} (which is equivalent to problem \eqref{base_saddle_point_problem_G}), we can try to use some gradient type algorithms, at each iteration of the used algorithm, we calculate the inexact gradient of $\phi(x, \cdot)$ for each $x \in \mathbb{R}^{d_x}$ in order to solve the (inner) maximization problem in \eqref{base_x_problem}. 

\rev{In practice, Multi-step Gradient Descent Ascent (MGDA) algorithm and its modifications are widely used to solve the problem \eqref{base_saddle_point_problem} (see e.g. \cite{Goodfellow,Wasserstein_GANs,Towards_Adversarial_Attacks}), as important algorithms for the considered type of saddle-point problems.}  In \cite{Solving_non_convex_min_max}, authors demonstrate the effectiveness of MGDA on a class of games in which one of the players satisfies the PL-condition and another player has a general non-convex structure. At the same time, \cite{Yang_vr_opt} shows that one step of gradient descent ascent demonstrates good performance and has theoretical guarantees for convergence in two-sided PL-games.

Finally, we mention that the main difference between distributed problems from usual optimization problems is to keep every agent's vectors to their average. It is approached through communication steps, where the communication can be performed in different scenarios. In our work, we will consider the time-varying network with changeable edges set and use the standard consensus procedure (see Subsec. \ref{consensus_procedure}) when the agent's vectors are averaged by multiplying by a weighted matrix of the graph at the current moment.

\subsection{Related works}
The decentralized algorithm makes two types of steps: local updates and information exchange. Local steps may use gradient \cite{nedic2017achieving,scaman2018optimal,pu2020push,qu2019accelerated,ye2023multi,li2018sharp} or sub-gradient \cite{nedic2009distributed} computations. In primal-only methods, the agents compute gradients of their local functions and alternate taking gradient steps and communication procedures. Under cheap communication costs, it may be beneficial to replace a single consensus iteration with a series of information exchange rounds. Such methods as MSDA \cite{scaman2017optimal}, D-NC \cite{jakovetic2018unification}, and Mudag \cite{ye2023multi} employ multi-step gossip procedures. In order to achieve acceptable complexity bounds, one may distribute accelerated methods directly \cite{qu2019accelerated,ye2023multi,li2018sharp,jakovetic2018unification,dvinskikh2021decentralized} or use a Catalyst framework \cite{li2020revisiting}. Accelerated methods meet the lower complexity bounds for decentralized optimization \cite{hendrikx2021optimal,li2020optimal,scaman2017optimal}. \rev{As usual, by complexity we mean a sufficient number of iterations of the algorithm that guarantee the solution of the problem with a given accuracy.}


Consensus restrictions $x_1 = \ldots = x_n$ may be treated as linear constraints, thus allowing for a dual reformulation of problem \eqref{main_sum_problem}. Dual-based methods include dual ascent and its accelerated variants \cite{scaman2017optimal,wu2019fenchel,zhang2017distributed,uribe2020dual}. Primal-dual approaches like ADMM \cite{arjevani2020ideal,wei2012distributed} are also implementable in decentralized scenarios. In \cite{scaman2018optimal}, the authors developed algorithms for non-smooth convex objectives and provided lower complexity bounds for this case, as well. 

Changing topology for time-varying networks requires new approaches to decentralized methods and a more complicated theoretical analysis. The first method with provable geometric convergence was proposed in \cite{nedic2017achieving}. Such primal algorithms as the Push-Pull Gradient Method \cite{pu2020push} and DIGing \cite{nedic2017achieving} are robust to network changes and have theoretical guarantees of convergence over time-varying graphs. Recently, a dual method for time-varying architectures was introduced in \cite{maros2018panda}. 

In \cite{rogozin_accel_do}, it was studied the problem of decentralized optimization with strongly convex smooth objective functions. The authors investigated accelerated deterministic algorithms under time-varying network constraints with a consensus projection procedure. The distributed stochastic optimization over time-varying graphs with consensus projection procedure was studied in \cite{rogozin2021accelerated}. The consensus projection procedure made it possible to use more acceptable parameters of smoothness and strong convexity in the complexity estimates. The main goal of this paper is to investigate a similar consensus projection approach for problems with PL-condition. Note that approach \cite{rogozin2021accelerated,rogozin_accel_do} is based on the well-known concept of the inexact $(\delta, L, \mu)$-oracle. But in the PL-case,  we will use inexact gradients with additive noise to describe the consensus projection procedure.

\subsection{Our contributions} \label{S:contrib}
Summing up, the contribution of this paper is as follows.
\begin{itemize}
	\item We study the sum-type minimization problem when the objective function satisfies the PL-condition. To solve a such class of problems, we propose a gradient descent type method with consensus projection procedure (see Algorithm \ref{alg:base_alg}) and access to only inexact gradient. We consider two cases of gradient inexactness: the bounded deterministic inexactness and the sum of random noise with deterministic bias. We estimated the sufficient communication steps and iterations number to approach the required quality concerning the function and the distance between the agent's vectors and their average for both cases.
    \item We study the decentralized saddle-point problem (with the structure of a sum) when the objective function satisfies the two-sided PL-condition. For solving a such generalized class of problems, we proposed a generalization of the  MGDA method (see Algorithm \ref{alg:base_alg_sp}) with a consensus procedure.  We provided an estimation for the sufficient number of iterations for inner and outer loops to approach the acceptable quality concerning the function. Also, we estimate the communication complexity to a distance between the agent's vectors and their average for both cases to be small enough. Additionally, we research the influence of interference on the convergence of the proposed Algorithm \ref{alg:base_alg_sp}. We suppose that the inexactness in the gradient can be separated into deterministic noise with the bounded norm and zero-mean noise with the finite second moment. 
    \item We present some numerical experiments, which demonstrate the effectiveness of the proposed algorithms, for the Least Squares and Robust Least Squares problems.
\end{itemize}

\section{Fundamentals and Assumptions for the problems under consideration}

Throughout the paper, $\langle \cdot, \cdot \rangle$ denotes the inner product of vectors or matrices. Correspondingly, by $\| \cdot \|$, we denote the $2$-norm for vectors or the Frobenius norm for matrices.


At the first, for the minimization problem \eqref{main_sum_problem} (or its equivalent \eqref{equiv_main_sum_problem}), we assume
\begin{assumption}[Lipschitz smoothness]\label{L_smooth_assump}
For every $i = 1, \ldots, n$, the function $f_i$ is $L_i$-smooth, for some $L_i >0$. 
\end{assumption}

Under this assumption, we find that the function $F(\mathbf{X})$ (see \eqref{function_FX}) is $L_l$-smooth on $\mathbb{R}^{n\times d}$, where $L_l=\max\limits_{1\leq i \leq n}  L_i$, and $L_g$-smooth on  $\mathcal{C}$, where $L_g =\frac{1}{n} \sum\limits_{i=1}^n L_i$. The constant $L_l$ is called a local constant and $L_g$ is called a global constant.

\begin{assumption}[PL-condition]\label{PL_assump}
The function $f$ satisfies the PL-condition, i.e., it holds the following inequality
\begin{equation}\label{PL_cond}
f(x)-f^* \leq \frac{1}{2\mu} \|\nabla f(x)\|^2, \quad \forall x \in \mathbb{R}^{d},
\end{equation}
for some $\mu >0$, and $f^*$ is the optimal value of the function $f$.
\end{assumption}

Also under this assumption, we find that $f$  satisfies the quadratic growth condition (QG-condition) \rev{(see \cite{karimi2016linear}}):
\begin{equation}\label{QC_cond}
\|x-x^*\|^2\leq \frac{2}{\mu}\left(f(x)-f^*\right); \quad \forall x \in \mathbb{R}^d,
\end{equation}
where $x^*$ is the nearest point to the optimal solution of the minimization problem under consideration.

At the second, for the saddle-point problem \eqref{base_saddle_point_problem}, let us introduce the following assumption.

\begin{assumption}[Lipschitz smoothness]\label{L_smooth_assump_sp}
For every $i = 1, \ldots, n$, the function $\phi_i(\cdot,\cdot)$ is differentiable with respect to its both variables and smooth, i.e., the following inequalities hold
\begin{equation*}\label{cond}
\begin{aligned}
\left\|\nabla_x \phi_i(x_1, y_1)-\nabla_x \phi_i(x_2, y_2)\right\|\leq L_{xx,i} \|x_1-x_2\| + L_{xy,i} \|y_1-y_2\|, \\ 
\left\|\nabla_y \phi_i(x_1, y_1)-\nabla_y \phi_i(x_2, y_2)\right\|\leq L_{yx,i} \|x_1-x_2\| + L_{yy,i} \|y_1-y_2\|,
\end{aligned} 
\end{equation*}
for all $x_1, x_2 \in \mathbb{R}^{d_x}, y_1, y_2 \in \mathbb{R}^{d_y}$ and $ L_{xx,i},  L_{yy,i},  L_{xy,i},  L_{yx,i} \in \mathbb{R}_+^*$.
\end{assumption}

Note, that if Assumption \ref{L_smooth_assump_sp} is satisfied, then it will be also satisfied for the function $\Phi(\cdot, \cdot)$ on  $\mathbb{R}^{n \times d_x}\times \mathbb{R}^{n \times d_y}$ with constants $L_{xx,g}, L_{yx,g}, L_{xy,g}, L_{yy,g}$,  and on  $\mathcal{C}_{x}\times \mathcal{C}_{y}$, with constants $L_{xx,l}, L_{yx,l}, L_{xy,l}, L_{yy,l}$, where $L_{ab, l} =\max\limits_{1 \leq i \leq n} L_{ab, i}$,  $L_{ab, g} =\frac{1}{n} \sum\limits_{i=1}^n L_{ab, i}$, and $ab \in \{xx, yy, xy, yx\}$.

\subsection{Inexact gradient oracle}
We will divide the inexactness of the gradient into random noise and a deterministic bias. To describe this noise, we will use the following definition.
\begin{definition}[$(\delta, \sigma^2)$-biased gradient oracle]\label{BGO_def}
Let $\delta>0, \sigma >0$ and $\xi$  be a random variable with probability distribution $\mathcal{D}$. $(\delta, \sigma^2)$-biased gradient oracle is a map $\mathbf{g}:\mathbb{R}^d\times \mathcal{D}\rightarrow \mathbb{R}^d$, such that
\begin{equation*}
\left\|\mathbb{E}_{\xi}\mathbf{g}(x, \xi)-\nabla f(x)\right\|\leq \delta,\quad     \mathbb{E}_{\xi}\left\|\mathbf{g}(x, \xi)-\mathbb{E}_{\xi}\mathbf{g}(x, \xi)\right\|^2\leq \sigma^2.
\end{equation*}
\end{definition}
In other words, we suppose that the inexactness of the gradient contains two parts: bias and noise. At that, bias has a bounded norm and noise has bounded the second moment.

\subsection{Consensus procedure}\label{consensus_procedure}
We consider a sequence of non-directed communication graphs $\left\{\mathcal{G}^k = \left(V, E^k\right)\right\}_{k=0}^{\infty}$ and a sequence of corresponding mixing matrices $\{\mathbf{W}^k\}_{k=0}^{\infty}$ associated with it. We assume the following assumptions.
\begin{assumption}\label{w_mixing_matrix}
Mixing matrix sequence $\{\mathbf{W}^k\}_{k=0}^{\infty}$ satisfies the following properties:
\begin{itemize}
    \item (Decentralized property) If $(i,j)\notin E^k$, then $[\mathbf{W}^k]_{ij}=0$.
		
    \item (Double stochasticity) $\mathbf{W}^k \mathbf{1}_n=\mathbf{1}_n$, and  $\mathbf{1}_n^\top \mathbf{W}^k = \mathbf{1}_n^\top$.
		
    \item (Contraction property) There exist $\tau\in\mathbb{Z}_{++}$ and $\lambda\in(0,1)$ such that for every $k\geq \tau-1$, it holds the following inequality 
    $$
    \left\| \mathbf{W}^k_{\tau} \mathbf{X} - \overline{\mathbf{X}} \right\|\leq (1-\lambda) \left\| \mathbf{X}  - \overline{\mathbf{X}} \right\|,
    $$
    where $\mathbf{W}^k_{\tau}=\mathbf{W}^k\dots \mathbf{W}^{k-\tau+1}$.
\end{itemize}
\end{assumption}
\rev{The last property in Assumption \eqref{w_mixing_matrix} is a generalization of several well-known cases: time-static connected graph, sequence of the connected graph and $\tau$-connected graph sequence \cite{nedic2017achieving}. A stochastic variant of this contraction property is also studied in \cite{koloskova2020unified}.}

During every communication round, the agents exchange information according to the rule
$$
x_{i}^{k+1}=w_{i i}^{k}+\sum_{(i, j) \in E^{k}} w_{i j}^{k} x_{j}^{k} \text {. }
$$
In matrix form, this update rule writes as $\mathbf{W}^{k+1} = \mathbf{W}^{k} \mathbf{X}^{k}$. The contraction property in Assumption \eqref{w_mixing_matrix} is needed to ensure geometric convergence of Algorithm \ref{alg:consensus} to the average of nodes’ initial vectors, i.e., to $\overline{x}_0$. In particular, the contraction property holds for $\tau$-connected graphs with Metropolis weights choice for $\mathbf{W}^{k}$, i.e., 
\begin{equation}\label{matrix_W_metropolis}
\left[\mathbf{W}^{k}\right]_{i j}= 
\begin{cases}
1 / \left(1+\max \{d_{i}^{k}, d_{j}^{k}\}\right) & \text{if}\;\; (i, j) \in E^{k},
\\ 0 & \text{if} \;\; (i, j) \notin E^{k}, 
\\ 1-\sum_{(i, m) \in E^{k}} \left[ \mathbf{W}^{k} \right]_{i m} & \text{if} \;\; i=j,
\end{cases}
\end{equation}
where $d_i^k$ denotes the degree of node $i$ in graph $\mathcal{G}^k$.

\begin{algorithm}[htp]
	\caption{Consensus.} \label{alg:consensus}
	\begin{algorithmic}[1]
		\REQUIRE Initial point $\mathbf{Z}^0 \in \mathcal{C}$,  number of communication rounds $T$.
        \STATE Take current time moment $t_0$ from global variable.
		\FOR{$k =0, \ldots, T-1$}
		\STATE $\mathbf{Z}^{k+1}:= \mathbf{W}^{t_0+k}\, \mathbf{Z}^k$. 
		\ENDFOR 
		\STATE Update global variable with current time moment: $t_0=t_0+T$.
		\RETURN $\mathbf{Z}^T$.
	\end{algorithmic}
\end{algorithm}

\rev{Note, that the contraction property guarantee that after $T=N\tau$ steps of Algorithm \ref{alg:consensus} one obtains the point $\mathbf{Z}^0$ such that $$\|\mathbf{Z}^T-\overline{\mathbf{Z}}^0\|=\|W^{N\tau+t_0}\dots W^{t_0} \mathbf{Z}^0 - \overline{\mathbf{Z}}^0\|\leq (1-\lambda)^N \|\mathbf{Z}^0 -\overline{\mathbf{Z}}^0\|.$$ In other words, the consensus procedure converges with any accuracy because of contraction property in Assumption \ref{w_mixing_matrix}.}

\section{Algorithm for PL-minimization problem}

In this section, we focus on the minimization problem \ref{equiv_main_sum_problem} (which is equivalent to problem \eqref{main_sum_problem}). For this, we propose an algorithm (listed as Algorithm \ref{alg:base_alg}) using the inexact gradient. The proposed algorithm is a gradient-type method, that uses a $(\delta,0)$-biased gradient oracle $\widetilde{\nabla}F(\mathbf{X})$ without noise. This condition can be rewritten as  $\left\|\widetilde{\nabla}F(\mathbf{X})-{\nabla}F(\mathbf{X})\right\|\leq \delta$, where $\nabla F(\mathbf{X}) = \left(\nabla f_1(x_1), \ldots, \nabla f_n(x_n)\right)^{\top} \in \mathbb{R}^{d \times n}$ denotes the gradient of $F$ at $\mathbf{X}$. Note, that the considered inexact gradient in Algorithm \ref{alg:base_alg} is a usual additively inexact gradient. \rev{Further, we note $\overline{\mathbf{X}}^k$ and $\overline{\widetilde{\nabla} F}(\mathbf{X}^k)$ as averaged $\mathbf{X}^k$ and $\widetilde{\nabla} F(\mathbf{X}^k)$ over consensus.}

\begin{algorithm}[htp]
	\caption{Decentralized gradient descent with consensus subroutine.}
	\label{alg:base_alg}
	\begin{algorithmic}[1]
		\REQUIRE Starting point $\mathbf{X}^0 \in \mathcal{C}$,  step size $\gamma>0$, number of steps $N$, the sequence of number for communication rounds $\{T^k\}_{k=0}^{N-1}$ in consensus.
		\FOR{$k=0, \ldots, N-1$}
		\STATE $\mathbf{Z}^{k+1}= \mathbf{X}^k - \gamma \widetilde{\nabla} F(\mathbf{X}^k)$,
		\STATE $\mathbf{X}^{k+1}= \text{Consensus}(\mathbf{Z}^{k+1}, T^k)$.
		\ENDFOR 
		\RETURN $\mathbf{X}^N$.
	\end{algorithmic}
\end{algorithm}

    Indeed, after the $k$-th iteration of Algorithm \ref{alg:base_alg} we have a point $\mathbf{X}^{k+1}\approx \overline{\mathbf{X}}^{k+1}=\overline{\mathbf{X}}^k-\gamma \overline{\widetilde{\nabla} F}(\mathbf{X}^k)$. If we assume that at each iteration, $\mathbf{X}^k$ is close enough to its projection on $\mathcal{C}$, i.e., $\left\|\mathbf{X}^k-\overline{\mathbf{X}}^k\right\|^2\leq \delta^{\prime}$ for some $\delta^{\prime}>0$, then we can estimate the new inexactness. So, in the fact, it is a variant of gradient method for problem \eqref{equiv_main_sum_problem} with an additively inexact gradient. At the same time, $\overline{x}^{k+1}=\overline{x}^k-\gamma \overline{\widetilde{\nabla} f}(x^k)$ is the usual gradient method for minimizing the function $f$ (i.e., for the problem \eqref{main_sum_problem}) with an additively inexact gradient $\overline{\widetilde{\nabla} f}(x^k)=\frac{1}{n}\sum\limits_{i=1}^n\widetilde{\nabla}f_i(x_{i}^k)$.

Let us introduce the following constant characterizes the distance from consensus space to points generated by Algorithm.
\begin{equation}\label{def_D}
\sqrt{D}=\gamma \left\|\nabla F(\overline{\mathbf{X}}^*) \right\|+\sqrt{\delta'} + \left(\gamma+\frac{1}{\mu} \right) \Delta + \gamma L_g \sqrt{\frac{2}{\mu} \left(1-\frac{\mu}{L_g}\right) \left(F(\overline{\mathbf{X}}^0)-F^*\right)}, 
\end{equation}
where  $\gamma=\frac{1}{L_g}$. Using this value we can prove the following theorem.

\begin{theorem}\label{theorem:main_theorem}
Choose some $\varepsilon >0$, define $D>0$ as in \eqref{def_D}, and set $\Delta = \delta+L_l\sqrt{\delta'}, \gamma=\frac{1}{L_g}$. Under Assumptions \eqref{L_smooth_assump} and \eqref{PL_assump}, Algorithm \ref{alg:base_alg} requires
$$
N = \left\lceil\frac{L_g}{\mu}\log\left(\frac{f(\overline{x}^0)-f^*}{\varepsilon}\right)\right\rceil
$$
gradient computation at each node, $T=\tau\left\lceil\frac{1}{2\lambda}\log \frac{D}{\delta'}\right\rceil$ communication steps at each iteration and $N_{tot}=N\cdot T$ communication steps to yield $\mathbf{X}^N$ such that
\begin{equation}
f(\overline{x}^N)-f^*\leq \varepsilon+\frac{\Delta^2}{2\mu n}, \quad \text{and} \quad \left\|\mathbf{X}^N-\overline{\mathbf{X}}^N \right\|\leq \sqrt{\delta^{\prime}}.
\end{equation}
\end{theorem}
\begin{proof}
    Firstly, let us estimate the inexactness for inexact gradient $\overline{\widetilde{\nabla}F}({\mathbf{X}}^k)$ in the following way:
$$
\left\|\overline{\widetilde{\nabla}F}({\mathbf{X}}^k)-\overline{{\nabla}F}(\overline{\mathbf{X}}^k)\right\|\leq \delta + L_l\sqrt{\delta^{\prime}}.
$$

Using this result we can obtain the similar estimate for $\overline{\widetilde{\nabla} f}(x^k)=\frac{1}{n}\sum\limits_{i=1}^n\widetilde{\nabla}f_i(x_{i}^k)$:
\begin{equation}
\label{inexact_grad_x_f}
\left\|\overline{\widetilde{\nabla} f}(x^k)-\nabla f(\overline{x}^k)\right\|^2 = \frac{1}{n}\left\|\overline{\widetilde{\nabla}F}(\mathbf{X}^k)-\overline{{\nabla}F}(\overline{\mathbf{X}}^k)\right\|^2 = \frac{\left(\delta + L_l\sqrt{\delta^{\prime}}\right)^2}{n}.
\end{equation}

Using Assumption \ref{PL_assump}, inequality \eqref{inexact_grad_x_f} and taking step size $\gamma=\frac{1}{L_g}$, we can estimate the  convergence rate of Algorithm \ref{alg:base_alg} with respect to the function $f$, as follows
\begin{equation}
\label{convergence_rate}
f(\overline{x}^{k})-f^*\leq \left(1-\frac{\mu}{L_g}\right)^{k} \left(f(\overline{x}^0)-f^*\right) + \frac{\Delta^2}{2\mu n},
\end{equation}
where $\Delta = \delta+L_l\sqrt{\delta'}$. Although at each iteration of Algorithm \ref{alg:base_alg} we have access only to $\mathbf{X}^k$. 

Now, let us find the sufficient communication step such that the following expression is true
\begin{equation}\label{cond_delta}
\left\|\mathbf{X}^{k}-\overline{\mathbf{X}}^{k}\right\|\leq \sqrt{\delta^{\prime}} \Longrightarrow \left\|\mathbf{X}^{k+1}-\overline{\mathbf{X}}^{k+1}\right\|\leq \sqrt{\delta^{\prime}}.
\end{equation}

By the contraction property (see Assumption \eqref{w_mixing_matrix}) and using that $\overline{\mathbf{Z}}^{k+1}=\overline{\mathbf{X}}^{k+1}$, we have
\begin{equation}\label{ineq_1414}
\left\|\mathbf{X}^{k+1}-\overline{\mathbf{X}}^{k+1}\right\|\leq (1-\lambda)^{\left\lfloor\frac{T_k}{\tau}\right\rfloor} \left\|\overline{\mathbf{X}}^{k+1} - \mathbf{Z}^{k+1}\right\|.
\end{equation}

Let us estimate the right-hand side of inequality \eqref{ineq_1414}.
\begin{equation}\label{asdfg}
\begin{aligned}
\left\|\overline{\mathbf{X}}^{k+1}-\mathbf{Z}^{k+1}\right\| &\leq \left\|\overline{\mathbf{X}}^{k}-\mathbf{Z}^{k+1}\right\| \leq \left\|\overline{\mathbf{X}}^{k}-\mathbf{X}^k\right\| +  \gamma \left\| \widetilde{\nabla}F(\mathbf{X}^k)\right\|
\\& \leq 
\sqrt{\delta^{\prime}} + \gamma \left\|\widetilde{\nabla}F(\mathbf{X}^k) - \nabla F(\overline{\mathbf{X}}^k)\right\| + \gamma\left\|\nabla F(\overline{\mathbf{X}}^k)- \nabla F(\overline{\mathbf{X}}^*)\right\| 
\\& \;\;\;\; 
+ \gamma \left\|\nabla F(\overline{\mathbf{X}}^*)\right\| 
\\& \leq
\sqrt{\delta^{\prime}}+\gamma \Delta + \gamma L_g \left\|\overline{\mathbf{X}}^k -\overline{\mathbf{X}}^*\right\| + \gamma \left\|\nabla F(\overline{\mathbf{X}}^*)\right\|.
\end{aligned}
\end{equation}
From quadratic growth condition \eqref{QC_cond}, we have $$
\left\|\overline{\mathbf{X}}^k-\overline{\mathbf{X}}^*\right\|^2 = n \left\|\overline{x}^k-\overline{x}^* \right\|^2 \leq \frac{2n}{\mu} \left( f( \overline{x}^k) - f^*\right).
$$

Further, using the convergence rate \eqref{convergence_rate} we can estimate the value of $\left\|\overline{\mathbf{X}}^k-\overline{\mathbf{X}}^*\right\|^2$ in \eqref{asdfg}, as the following
$$
\left\| \overline{\mathbf{X} }^k - \overline{\mathbf{X}}^* \right\|^2 \leq \frac{2}{\mu} \left(1 - \frac{\mu}{L_g} \right)^{k+1} \left(F(\overline{\mathbf{X}}^0) - F^*\right) + \frac{\Delta^2}{\mu^2}.
$$
Therefore we have $\left\| \overline{\mathbf{X}}^{k+1} - \mathbf{Z}^{k+1} \right\| \leq \sqrt{D}$, where $D$ is the constant defined according to \ref{def_D}

Thus, for $T_k=\tau\left\lceil\frac{1}{2\lambda}\log \frac{D}{\delta'}\right\rceil$, \eqref{cond_delta} will be satisfied. So, uniting \eqref{convergence_rate} and the result about communication steps per iteration gets the theorem statement.
\end{proof}

Further, we consider the case when Algorithm \ref{alg:base_alg} uses $(\delta, \sigma^2)$-biased gradient oracle $\widetilde{\nabla}F(\mathbf{X})$ for arbitrary values $\delta$ and $\sigma$.

As in the previous, after the $k$-th iteration of Algorithm \ref{alg:base_alg} we have a point $\mathbf{X}^{k+1} \approx \overline{\mathbf{X}}^{k+1} = \overline{\mathbf{X}}^k - \gamma \overline{\widetilde{\nabla} F} (\mathbf{X}^k)$. In this case, if we consider $\overline{\widetilde{\nabla} F}(\mathbf{X}^k)$ as an approximation of the exact gradient $\overline{\nabla F}(\overline{\mathbf{X}}^k)$, then will have three inexactness: inexactness caused by $\overline{\mathbf{X}}^k\neq \mathbf{X}^k$, bias at point $\mathbf{X}^k$ and zero-mean noise 
at $\mathbf{X}^k$. Note, that $\mathbf{X}^k$ is also a random variable.

Let us estimate the bias in $\overline{\widetilde{\nabla} F}(\mathbf{X}^k)$ at the given point $\overline{\mathbf{X}}^k$, as follows 
\begin{equation*}
\begin{aligned}
& \sqrt{n}\left\|\mathbb{E}_{x^k,\xi}\overline{\widetilde{\nabla} f}(x^k) -{{\nabla} f}(\overline{x}^k)\right\| 
= \left\|\mathbb{E}_{\mathbf{X}^k,\xi}\overline{\widetilde{\nabla} F}(\mathbf{X}^k)-\overline{{\nabla} F}(\overline{\mathbf{X}}^k)\right\|
\\& \qquad \qquad \qquad 
\leq \left\|\mathbb{E}_{\mathbf{X}^k, \xi}\overline{\widetilde{\nabla} F}(\mathbf{X}^k) - \mathbb{E}_{\mathbf{X}^k} \overline{{\nabla} F}(\mathbf{X}^k)\right\| + \left\|\mathbb{E}_{\mathbf{X}^k}\overline{{\nabla} F}(\mathbf{X}^k)-\overline{{\nabla} F}(\overline{\mathbf{X}}^k)\right\|
\\& \qquad \qquad \qquad
\leq \mathbb{E}_{\mathbf{X}^k}\left\|\mathbb{E}_{\xi}\overline{\widetilde{\nabla} F}(\mathbf{X}^k)-\overline{{\nabla} F}(\mathbf{X}^k)\right\| + \mathbb{E}_{\mathbf{X}^k}\left\|\overline{{\nabla} F}(\mathbf{X}^k)-\overline{{\nabla} F} (\overline{\mathbf{X}}^k)\right\|
\\& \qquad \qquad \qquad
\leq \delta + L_l\mathbb{E}_{\mathbf{X}^k} \left\|\mathbf{X}^k-\overline{\mathbf{X}}^k\right\|,
\end{aligned}
\end{equation*}
\rev{where $\mathbb{E}_{\mathbf{X}^k}$ and $\mathbb{E}_\xi$ mean conditional mathematical expectations under variable $\mathbf{X}^k$ and random variable $\xi$ for the given $\overline{\mathbf{X}}^k$. As a result, one requires $\mathbb{E}_{\mathbf{X}^k} \left\|\mathbf{X}^k-\overline{\mathbf{X}}^k\right\|^2\leq {\delta'}$ for small bias of gradient.} On the other hand, we can construct Algorithm \ref{alg:base_alg} in such a way that $\mathbb{E} \left\|\mathbf{X}^k-\overline{\mathbf{X}}^k\right\|^2 \leq {\delta'}$ at all iterations. When we assume, that Consensus procedure \ref{alg:consensus} guarantees  $\mathbb{E} \left\|\mathbf{X}^k-\overline{\mathbf{X}}^k\right\|^2 \leq {\delta'}$, we have the following estimation
$$
\mathbb{E}\left\|\mathbb{E}_{x,\xi} \overline{\widetilde{\nabla} f}(x^k)-{{\nabla} f}(\overline{x}^k)\right\|^2\leq \frac{2\delta^2 + 2L_l^2 \delta^{\prime}}{n}.
$$

At the same time, we can estimate the random noise in a similar way (see Appendix \ref{app:noise_est}) and it gets the following lemma.

\begin{lemma}\label{est_noise_lemma}
Let us assume that the functions $f_i$ meet Assumption \ref{L_smooth_assump} and for all $j\leq k$ we have that $\mathbb{E}\left\|\mathbf{X}^j-\overline{\mathbf{X}}^j\right\|^2\leq \delta^{\prime}$, where $\{\mathbf{X}^j\}_{j\leq k}$. Then for the bias and noise in the inexact gradient $\overline{\widetilde{\nabla} f}(x^k)$,  
we have the following estimations
\begin{equation*}
\mathbb{E}\left\|\mathbb{E}_{x,\xi} \overline{\widetilde{\nabla} f}(x^k)-{{\nabla} f}(\overline{x}^k)\right\|^2\leq \frac{2\delta^2 + 2L_l^2 \delta'}{n},
\end{equation*}
and
\begin{align*}
\mathbb{E}\left\|\mathbb{E}_{x,\xi}\overline{\widetilde{\nabla} f}(x^k)-\overline{\widetilde{\nabla} f}(x^k)\right\|^2\leq \frac{16L_l^2{\delta'}+18\sigma^2+16\delta^2}{n}.
\end{align*}
\end{lemma}

{\color{black}
Note that the mathematical expectation $\mathbb{E}$, above is not conditional. So, $\overline{\widetilde{\nabla} f}(x^k)$ is not a biased in the sense of Definition \ref{BGO_def}. Nevertheless, we can use the following gradient  method
\begin{equation}
\label{method}
x_k = x_{k-1}-\gamma g(x_k, \xi), \quad k = 1,2, \ldots 
\end{equation}
where $g(x_k, \xi)=\nabla f(x_k) + n(x_k,\xi)+b(x^k)$ such that $\mathbb{E}\left\|n(x_k,\xi)\right\|^2 \leq \sigma^2$ and $\mathbb{E}\left\|b(x_k)\right\|^2 \leq \delta^2$. The convergence of such a method is given by the following lemma (see Lemma 2 in \cite{stich_sgd_biased}).
}

\begin{lemma}\label{lemma_conv}
Let  $f$ be a function that satisfies the PL-condition for a constant $\mu > 0 $ and be an $L$-smooth function. The gradient oracle $g(x_k, \xi)=\nabla f(x_k) + n(x_k,\xi)+b(x^k)$ is such that $\mathbb{E}\left\|n(x_k,\xi)\right\|^2 \leq \sigma^2$ and $\mathbb{E}\left\|b(x_k)\right\|^2 \leq \delta^2$. When $\gamma\leq \frac{1}{L}$, we can guarantee that method \eqref{method} converges to $f^*$ in the following way
\begin{equation}
\label{convergence_rate_lemma}
\mathbb{E}\left[f(x_{k})-f^*\right]\leq \left(1-\gamma\mu\right)^{k} (f(x_0)-f^*) + \frac{\delta^2}{2\mu} + \frac{L\gamma \sigma^2}{2\mu}.
\end{equation}
\end{lemma}

In the similar way, like for proof of Theorem \ref{theorem:main_theorem}, we can estimate required consensus steps $T$ for condition $\mathbb{E}\left\|\mathbf{X}^k -\overline{\mathbf{X}}^k\right\|^2\leq \delta^{\prime}$ (see Appendix \ref{app:proof_main_theorem_full}). So, we can state the following theorem.

\begin{theorem}\label{theorem:main_theorem_full}
Let  $f$ be a function meets Assumption \ref{PL_assump}, the functions $f_i$ meet Assumption \ref{L_smooth_assump} and the map $\widetilde{\nabla}F$ be a $(\delta,\sigma^2)$-biased gradient oracle. Define value $D$ as
$$
D=6\gamma^2 \left\|\nabla F(\overline{\mathbf{X}}^*)\right\|^2 + 2\delta^{\prime} + 6\left(\gamma^2+\frac{1}{\mu^2}\right) \Delta^2 +  \frac{12 \gamma^2 L_g^2}{\mu} \left(1-\frac{\mu}{L_g}\right) \left(F(\overline{\mathbf{X}}^0)-F^*\right),
$$
where $\Delta^2=18\left(L_l^2{\delta'}+\sigma^2+\delta^2\right)$ and step-size $\gamma=\frac{1}{L_g}$. Then Algorithm \ref{alg:base_alg} requires $N = \left\lceil\frac{L_g}{\mu}\log\left(\frac{f(\overline{x}^0)-f^*}{\varepsilon}\right)\right\rceil$ gradient computation at each node, $T=\tau\left\lceil\frac{1}{2\lambda}\log \frac{D}{\delta'}\right\rceil$ communication steps in each iteration and $N_{tot}=N \cdot T$ communication steps to yield $\mathbf{X}^N$ such that
\begin{equation}\label{quality_x}
\mathbb{E}f(\overline{x}^N)-f^*\leq \varepsilon+\frac{\Delta^2}{2\mu n}, \quad \text{and} \quad  \mathbb{E}\left\|\mathbf{X}^N - \overline{\mathbf{X}}^N\right\|^2 \leq {\delta^{\prime}}.
\end{equation}
\end{theorem}

\begin{remark}
By using the step-size $\gamma < \frac{1}{L_g}$   small enough, we can improve the accuracy level on function in \eqref{quality_x}. According to Lemma \ref{lemma_conv}, for step-size $\gamma$, Algorithm \ref{alg:base_alg} can converge with the rate as in the inequality \eqref{quality_x}, where 
$$
\Delta^2 = 2\delta^2 + L_l^2 \delta^{\prime} + L_g \gamma \left(16L_l^2{\delta^{\prime}}+18\sigma^2+16\delta^2\right).
$$
But in this case the number of required iterations $N$ increases in $\frac{L_g}{\gamma}$ times. Note, that parameter $D$ does not change.
\end{remark}

\begin{remark}
    Note, that the convergence rate is almost optimal for functions with PL-condition (see \cite{PL_optimality}).
\end{remark}

\section{Algorithm for distributed saddle point problems: PL--PL case}

In this section, we will consider a generalization of Algorithm \ref{alg:base_alg} (see Algorithm \ref{alg:base_alg_sp}) for the saddle-point problem \eqref{base_saddle_point_problem} (or its equivalent \eqref{base_saddle_point_problem_G}).

Additionally, we will research the influence of {\color{black} the inexact access to the oracle (especially using the inexact information of the gradient)} on the convergence of the proposed Algorithm \ref{alg:base_alg_sp}. We suppose, that the inexactness in the gradient can be separated into deterministic noise with the bounded norm and zero-mean noise with the finite second moment. Note, that the influence of such inexactness for Gradient Descent Ascent is well researched both in convex-concave and in non-convex-non-concave saddle point problems under various assumptions (see e.g. \cite{Yang_vr_opt}, \cite{Beznosikov_sgda_unified_theory}). Also, the results of \cite{Solving_non_convex_min_max} can be generalized for stochastic oracle (see Remark 3.8 in \cite{Solving_non_convex_min_max}).

At each iteration, the proposed method optimizes by an inner variable $\mathbf{Y}$ for the fixed outer variable $\mathbf{X}$ by Algorithm \ref{alg:base_alg}. After that, it makes one step of Algorithm \ref{alg:base_alg_sp} by outer variable. So, we obtain the well-known Multi-step Gradient Descent Ascent method with a consensus subroutine after each gradient method step. 

\begin{algorithm}[htp]
	\caption{Multi-step Gradient Descent Ascent with consensus (MGDA).}
	\label{alg:base_alg_sp}
	\begin{algorithmic}[1]
		\REQUIRE Number of the outer and inner steps $N_x, N_y$, starting point $(\mathbf{X}^0, \mathbf{Y}^0)$,  step sizes $\gamma_x$ and $\gamma_y$, the sequences of steps $\{T^k_x\}_k$ and $\{T^{k,j}_y\}_{k,j}$ in consensus.
		\FOR{$k =0,\ldots, N_x-1$}
		\STATE $\hat{\mathbf{Y}}^0 := \mathbf{Y}^k$.
		\FOR{$j=0,\ldots, N_y-1$}
		\STATE $\mathbf{Z}_y^{j+1}:= {\mathbf{\hat{Y}}}^j + \gamma_y \widetilde{\nabla}_{\mathbf{Y}} \Phi(\mathbf{X}^k, \hat{\mathbf{Y}}^j)$.
		\STATE $\hat{\mathbf{Y}}^{j+1}:= \text{Consensus}(\mathbf{Z}_y^{j+1}, T^{k,j}_y)$.
		\ENDFOR
		\STATE ${\mathbf{Y}}^{k+1} := \hat{\mathbf{Y}}^{N_y}$.
		\STATE ${\mathbf{Z}}_x^{k+1}:= \mathbf{X}^k - \gamma_x \widetilde{\nabla}_{\mathbf{X}} \Phi(\mathbf{X}^k, \mathbf{Y}^{k+1})$.
		\STATE ${\mathbf{X}}^{k+1}:= \text{Consensus}(\mathbf{Z}_x^{k+1}, T^k_x)$. 
		\ENDFOR 
		\RETURN $\mathbf{X}^{N_x}, {\mathbf{Y}}^{N_x}$.
	\end{algorithmic}
\end{algorithm}

We suppose that functions $\phi(\cdot, y)$ and $-\phi(x, \cdot)$ satisfy the PL-condition (Assumption \ref{PL_assump}). So, for each $x$ and $y$ the following inequalities hold
\begin{align}
\label{pl_cond_f}
\phi(x, y) - \min_{x \in \mathbb{R}^{d_x}}  \phi(x, y) \leq \frac{2}{\mu_x}\|\nabla_x \phi(x, y)\|^2,\\
\label{pl_cond_g}
\max_{y \in \mathbb{R}^{d_y}} \phi(x, y) -  \phi(x, y) \leq \frac{2}{\mu_y}\|\nabla_y \phi(x, y)\|^2.
\end{align}

Now, let us introduce functions in $\mathbf{Y}$ for fixed variable $\mathbf{X}$: $g_{\overline{x}}(y)= \phi(\overline{x}, y)$ and   $G_{\overline{\mathbf{X}}}(\mathbf{Y})= \Phi(\overline{\mathbf{X}}, \mathbf{Y})$. For the maximization of $g_{\overline{x}}(y)$ we define  $g^*_x := \max_{y\in\mathbb{R}^{d_y}} g_{\overline{x}}(y)$.

Additionally, we assume that there are points for functions $F$ and $G_{\overline{X}}$ in consensus subspace minimizing this functions in the full space. In other words, the following conditions hold:

\begin{equation}
    \label{y_minimizer}
    \forall \overline{X}\, \exists \overline{Y}^* \in \mathcal{C}_y: \max_{Y\in \mathbb{R}^{n \times d_y}} G_{\overline{X}}(Y) = G_{\overline{X}}(\overline{Y}^*),
\end{equation}
and
\begin{equation}
    \label{x_minimizer}
    \exists \overline{X}^* \in \mathcal{C}_x: \min_{X\in \mathbb{R}^{n \times d_x}} F(X) = F(\overline{X}^*).
\end{equation}
This conditions allows to obtain that Assumption \ref{PL_assump} holds for full space for any subproblem.

Note, in the inner iterations we have gradient process with respect to $\mathbf{Y}$. Uniting this assumption and results of the previous part we obtain the following result.

\begin{lemma}
\label{Y_lemma_convergence}
Let us define $\Delta_y^2 = \delta +L_{yy,l} \sqrt{\delta'_y} + L_{yx,l} \sqrt{\delta'_x} $ and $D_y$ by the following way
\begin{equation}\label{D_x^y_est}
\begin{aligned} 
\sqrt{D_{\mathbf{X},\mathbf{Y}}}&=\gamma_y{\left\|\nabla G_{\overline{x}} (\overline{\mathbf{Y}}^*) \right\|} + \sqrt{\delta'_y}+ \left(\gamma_y+\frac{1}{\mu_y}\right) \Delta_y 
\\& \;\;\;\; + \gamma_y{L_{yy,g}} \sqrt{\frac{2n}{\mu_y} \left(1-\frac{\mu_{y}}{L_{yy,g}}\right) \left(G_{\overline{x}}(\overline{\mathbf{Y}})-G_{\overline{x}}^*\right)}.
\end{aligned}
\end{equation}
Besides let us method makes at least $N_y = \left\lceil\frac{L_{yy,g}}{\mu}\log\frac{g_x (\overline{y}^0)-g_x^*}{\varepsilon}\right\rceil$ iterations of inner loop for $T_y=\tau\left\lceil\frac{1}{2\lambda}\log \left(\frac{D_{y,x}}{\delta'_y}\right)\right\rceil$ communication steps. If Assumption \ref{PL_assump} and statement \eqref{y_minimizer}
hold, the method obtains a point $\mathbf{Y}^k$ for any outer iteration $k$ such that    
\begin{equation*}
\max_{y \in \mathbb{R}^{d_y}} \phi(\overline{x}^k, \overline{y}^k) - \phi(\overline{x}^k, \overline{y}^k)\leq \varepsilon_y+\frac{\Delta_y^2}{2\mu_y n}=\hat{\varepsilon}_y, \quad \text{and} \quad \left\|\mathbf{Y}^k-\overline{\mathbf{Y}}^N\right\|\leq \sqrt{\delta'_y},
\end{equation*}
\end{lemma}

\begin{proof}
    Note, that in the inner iterations, we have the gradient process in the form 
$$
\hat{\mathbf{Y}}^j\approx \overline{\hat{\mathbf{Y}}}^j = \overline{\hat{\mathbf{Y}}}^{k-1} +\gamma_y \widetilde{\nabla}_{\mathbf{Y}} \Phi(\mathbf{X}^k, \hat{\mathbf{Y}}^j).
$$ 
Let $\widetilde{\nabla}_{\mathbf{Y}} \Phi(\mathbf{X}, \mathbf{Y})$ be a $(\delta_y, 0)$-biased gradient oracle for any fixed $\mathbf{X}$. So, we can decompose the bias component in $\widetilde{\nabla}_{\mathbf{Y}} \Phi(\mathbf{X}^k, \hat{\mathbf{Y}}^j)$ in the following way
\begin{equation}\label{sdfghh}
\begin{aligned}
\left\|{\widetilde{\nabla}_{\mathbf{Y}} \Phi}(\mathbf{X}^k, \hat{\mathbf{Y}}^j) - \nabla_y \Phi(\overline{\mathbf{X}}^k, \overline{\hat{\mathbf{Y}}}^j)\right\| \leq & \left\|{\widetilde{\nabla}_{\mathbf{Y}} \Phi}(\mathbf{X}^k, \hat{\mathbf{Y}}^j) - \nabla_y \Phi({\mathbf{X}}^k, {\hat{\mathbf{Y}}}^j)\right\|
\\&
+\left\|\nabla_y \Phi({\mathbf{X}}^k, {\hat{\mathbf{Y}}}^j) - \nabla_y \Phi({\mathbf{X}}^k, \overline{\hat{\mathbf{Y}}}^j)\right\|
\\&
+\left \|\nabla_y \Phi({\mathbf{X}}^k, \overline{\hat{\mathbf{Y}}}^j) - \nabla_y \Phi(\overline{\mathbf{X}}^k, \overline{\hat{\mathbf{Y}}}^j)\right\|.
\end{aligned}
\end{equation}

The first term on the right-hand side of inequality \eqref{sdfghh} is not more than $\delta$. Assuming that $\left\|\overline{\hat{\mathbf{Y}}^j}-\hat{\mathbf{Y}}^j \right\|^2 \leq \delta'_y$ and $\left\|\overline{\mathbf{X}}^k - \mathbf{X}^k\right\|^2 \leq \delta'_x $ for any $k,j$, we can estimate the last two terms in \eqref{sdfghh} as $L_{yy}\sqrt{\delta'_y} + L_{yx}\sqrt{\delta'_x}$. So, we have the following estimate for bias in the inexact gradient $\overline{\widetilde{\nabla}_{\mathbf{Y}} \phi}(\mathbf{X}^k, \hat{\mathbf{Y}}^j) = \frac{1}{n} \sum\limits_{i=1}^n\phi(x_{k,i}, \hat{y}^j_i)$:
\begin{align*}
\left\|\overline{\widetilde{\nabla}_{\mathbf{Y}} \phi}(\mathbf{X}^k, \hat{\mathbf{Y}}^j) - \nabla_y \phi(\overline{x}_k, \overline{\hat{y}}^j)\right\|
&\leq \frac{1}{n} \left\|{\widetilde{\nabla}_Y \Phi}(\mathbf{X}^k, \hat{\mathbf{Y}}^j) - \nabla_y \Phi(\overline{\mathbf{X}}^k, \overline{\hat{\mathbf{Y}}}^j)\right\|
\\ &
\leq \frac{(\delta+L_{yy,l}\sqrt{\delta'_y} + L_{yx,l}\sqrt{\delta'_x})^2}{n}.
\end{align*}
Finally, note, that $g_x$ is $L_{yy,g}$ smooth function. So, the convergence by $\mathbf{Y}$ is described by Theorem \ref{theorem:main_theorem_full}. So, we have the result of this lemma. 
\end{proof}

Note, in not decentralized case we have that the inexact gradient is close enough to exact one: $\overline{\nabla_{\mathbf{X}} \Phi}(\overline{\mathbf{X}}^k, \mathbf{Y}^k) \approx \overline{\nabla_{\mathbf{X}} F}(\overline{\mathbf{X}}^k)=\overline{\nabla_{\mathbf{X}} \Phi}(\overline{\mathbf{X}}^k, \overline{\mathbf{Y}}^*)$  (see Lemma A.6 from \cite{Solving_non_convex_min_max}). This case corresponds to every full-connected graph. Let us generalize this result for the decentralized problem statement.

\begin{theorem}\label{theorem:main_theorem_sp}
Let functions $\phi(\cdot, y)$ and $-\phi(x, \cdot)$ meet Assumption \ref{PL_assump} with constants $\mu_x >0$ and $\mu_y>0$, and functions $\phi_i$ meet the Assumption \ref{L_smooth_assump_sp}. Besides statements \eqref{y_minimizer} and \eqref{x_minimizer} hold. The gradient oracle $(\widetilde{\nabla}_{\mathbf{X}} \Phi(\mathbf{X},\mathbf{Y}), \widetilde{\nabla}_{\mathbf{Y}} \Phi(\mathbf{X},\mathbf{Y}))$ be $(\delta,0)$-biased gradient oracle. Also, Assumption \ref{w_mixing_matrix} holds. Let us introduce inexactness values $\Delta_y:= \delta+L_{yy,l}\sqrt{\delta'_y}+L_{yx,l}\sqrt{\delta'_x}$ and $\Delta_x=\delta +L_{xx,l}\sqrt{\delta'_x} + L_{xy, l}\left(\sqrt{\frac{{\varepsilon}_y}{2\mu_y}}+\frac{\Delta_y}{2\mu_y\sqrt{n}}+\sqrt{\delta'_y}\right).$ Also define values $D_{\mathbf{X}}$ according to \eqref{1hhg0} and $D_{\mathbf{Y}}=\max_k {D_{\mathbf{X}^k,\mathbf{Y}}}$, where $D_{\mathbf{X}^k,\mathbf{Y}}$ defined in  \eqref{D_x^y_est}.
	
Let Algorithm \ref{alg:base_alg_sp} makes $T_x=\tau\left\lceil\frac{1}{2\lambda}\log \left(\frac{D_X}{\delta'_x}\right)\right\rceil$ communication steps in each outer iteration and $T_y = \tau\left\lceil \frac{1}{2\lambda}\log \left(\frac{D_Y}{\delta'_y} \right)\right\rceil$ communication steps at each inner iteration, with step-sizes $\gamma_x=\frac{1}{L_x}$ and $\gamma_y=\frac{1}{L_{yy,g}}$. Then Algorithm \ref{alg:base_alg_sp} requires 
$$
N_x= \left\lceil\frac{L_x}{\mu_x}\log \left( \frac{F(\overline{\mathbf{X}}^0) - F^*}{\varepsilon_x} \right) \right\rceil
$$ 
outer iterations, where $L_x = L_{xx,g} + \frac{L_{xy,g}}{\mu_y}$, and 
$$
N_y = \left\lceil \frac{L_{yy,g}}{\mu_y} \log \left(\frac{G_{\overline{x}}^* - G_{\overline{x}}(\overline{\mathbf{X}}^0)}{\varepsilon_y}\right) \right\rceil
$$
inner iterations for each outer iteration to yield the pair $\left(\mathbf{X}^{N_x}, \mathbf{Y}^{N_y}\right)$, such that
\begin{equation}\label{quality_sp1}
f(\overline{x}^{N_x}) - f^*\leq \varepsilon_x+\frac{\Delta_x^2}{2\mu_x n},\quad \left\|\mathbf{X}^{N_x}-\overline{\mathbf{X}}^{N_x}\right\| \leq \sqrt{\delta'_x},
\end{equation}
and 
\begin{equation}\label{quality_sp2}
\max_{ y \in \mathbb{R}^{d_y}} \phi(\overline{x}^{N_x}, y)-\phi(\overline{x}^{N_x}, \overline{y}^{N_y})\leq \varepsilon_y+\frac{\Delta_y^2}{2\mu_y n}, \quad \left\|\mathbf{Y}^{N_y} - \overline{\mathbf{Y}}^{N_y} \right\| \leq \sqrt{\delta'_y}.
\end{equation}
\end{theorem}
\begin{proof}
    We proved, that after $N_y$ iterations of inner loop, we have such point $Y=\hat{\mathbf{Y}}^{N_y}$ that $g_x^*-g_x(\overline{y})\leq \varepsilon_y$ and $\left\|\mathbf{Y}-\overline{\mathbf{Y}}\right\|\leq \sqrt{\delta'_y}$. Using QG-condition for the function $g_x$, we can estimate the distance from obtained point to the optimal one as $\left\|\mathbf{Y}-\overline{\mathbf{Y}}^*(\mathbf{X})\right\| \leq \sqrt{\frac{\hat{\varepsilon}_y}{2\mu_y}} + \sqrt{\delta'_y}$,
where $\overline{y}^*(x)=\arg\max_{y\in\mathbb{R}^{d_y}} g_{\overline{x}}.$ So, we can estimate the inexactness of the gradient with respect to $\mathbf{X}$ at the point $\mathbf{Y} = \hat{\mathbf{Y}}^{N_y}$  according to the Assumption \ref{L_smooth_assump_sp} in the following way
\begin{align*}
\left\|\overline{\widetilde{\nabla}_{\mathbf{X}} \Phi}(\mathbf{X},\mathbf{Y})-\overline{\nabla F}(\overline{\mathbf{X}})\right\|
& 
\leq \delta+ L_{xx,l}\left\|\mathbf{X}-\overline{\mathbf{X}}\right\| + L_{xy,l} \left\|\mathbf{Y} - \overline{\mathbf{Y}}^*(\overline{\mathbf{X}}) \right\|
\\&
\leq \delta+L_{xx,l}\sqrt{\delta'_x} + L_{xy, l}\left(\sqrt{\frac{\hat{\varepsilon}_y}{2\mu_y}}+\sqrt{\delta'_y}\right)
\\&
\leq \delta +L_{xx,l}\sqrt{\delta'_x} + L_{xy, l}\left(\sqrt{\frac{{\varepsilon}_y}{2\mu_y}}+\frac{\Delta_y}{2\mu_y\sqrt{n}}+\sqrt{\delta'_y}\right).
\end{align*}

So, through the inner iterations, we obtain an inexact gradient with respect to $x$ such that $\left\|\widetilde{\nabla} f(x) - \nabla f(\overline{x})\right\|\leq \frac{\Delta_x}{\sqrt{n}}$, where
$$
    \Delta_x=\delta +L_{xx,l}\sqrt{\delta'_x} + L_{xy, l}\left(\sqrt{\frac{{\varepsilon}_y}{2\mu_y}}+\frac{\Delta_y}{2\mu_y\sqrt{n}}+\sqrt{\delta'_y}\right).
$$

At the same time, the function $\max_{y\in\mathbb{R}^{d_y}} \phi(x,y)$ meets PL-condition. Moreover, as known (see Lemma A.3 in \cite{Yang_vr_opt}), when the objective function $\phi$ meets Assumptions \ref{L_smooth_assump_sp} and condition \eqref{pl_cond_f}, the function $\max_{y\in\mathbb{R}^{d_y}} \phi(x,y)$ meets PL-condition with constant $\mu_x$ until if $\min_x \max_{y\in\mathbb{R}^{d_y}} \phi(x,y) = \min_{\mathbf{X}} F(\mathbf{X}).$. Also note, that under Assumptions \ref{L_smooth_assump_sp} and PL-condition \eqref{pl_cond_g} for inner variable, the function $F(\mathbf{X})=\Phi(\mathbf{X}, \mathbf{Y}^*(\mathbf{X}))$ is an $L_x$-Lipschitz continuous function with $L_x=L_{xx,g} + \frac{L_{xy,g}}{\mu_y}$ (see Lemma A.2 in \cite{Yang_vr_opt}). So, from Theorem \ref{theorem:main_theorem}, we obtain that Algorithm \ref{alg:base_alg_sp} requires $N_x = \left\lceil\frac{L_x}{\mu_x} \log \left(\frac{F(\overline{\mathbf{X}}^0)-F^*}{\varepsilon_x}\right) \right\rceil$ outer iterations and $T_x=\tau\left\lceil\frac{1}{2\lambda}\log \left(\frac{D_{\mathbf{X}}}{\delta'_x}\right)\right\rceil$ where 
\begin{equation}\label{1hhg0}
\sqrt{D_{\mathbf{X}}}  = \gamma_x{\left\|\nabla F(\overline{\mathbf{X}}^*)\right\|}+\sqrt{\delta'_x} + \left(\gamma_x+\frac{1}{\mu_x}\right) \Delta_x + \gamma_x L_x\sqrt{\frac{2}{\mu_x}\left(1-\frac{\mu_{x}}{L_{x}}\right) \left(F(\overline{\mathbf{X}}^0)-F^*\right)}.
\end{equation}

Uniting obtained above results for convergences by $\mathbf{X}$ and $\mathbf{Y}$, we obtain the result of the theorem.
\end{proof}

\begin{remark}
In general, Algorithm \ref{alg:base_alg_sp} requires
$$
N_x \cdot N_y = O \left( \left( \frac{L_{xx,g}L_{yy,g}}{\mu_y\mu_x}+\frac{L_{xy,g} L_{yy,g}} {\mu_y^2\mu_x} \right)\log^2 \frac{1}{\varepsilon} \right),
$$
gradient computations with respect to $\mathbf{Y}$ and
$$
N_x=O\left(\left(\frac{L_{xx,g}}{\mu_x}+\frac{L_{xy,g}}{\mu_y\mu_x}\right)\log\frac{1}{\varepsilon}\right),
$$ 
gradient computations with respect to $\mathbf{X}$ at each node.
\end{remark}

\begin{remark}
Algorithm \ref{alg:base_alg_sp} requires $T_{tot}=N_xT_x + N_yN_xT_y$ communication steps to achieve the required quality. Note, that the communication steps in inner iterations and outer iterations can have different computational costs. Namely, the dimensions $d_x$ for $\mathbf{X}$ and $d_y$ for $\mathbf{Y}$ can significantly differ.
\end{remark}

\begin{remark}
\rev{    The main restriction in the obtained result is conditions \ref{y_minimizer} and \ref{x_minimizer} and PL-condition. Nevertheless, note that such situation is typical for overparametrized problems that includes many problems from modern deep learning (see \cite{fit_without_fear}, \cite{karimi2016linear}).}
\end{remark}

\begin{remark}
    \rev{The particular case of functions with PL-condition is strong convexity. In such case, the result can be significantly accelerated (see e.g. \cite{rogozin_accel_do}).}
\end{remark}

Now, we consider the case, when $\widetilde{\nabla}_{\mathbf{X}} \Phi(\mathbf{X},\mathbf{Y})$ and $\widetilde{\nabla}_{\mathbf{Y}} \Phi(\mathbf{X},\mathbf{Y})$ be $(\delta,\sigma^2)$-biased gradient oracle. Let parameters numbers of communication steps $T_x$ and $T_y$ such that $\mathbb{E}\|\mathbf{X}^j-\overline{\mathbf{Y}}^j\|^2 \leq {\delta'_x}$ and $\mathbb{E}\|\mathbf{Y}^j-\overline{\mathbf{Y}}^j\|^2 \leq {\delta'_y}$ in all iterations. So, we can estimate noise and bias for the gradient with respect to $\mathbf{Y}$ for each inner loop. The proof of the following lemma is presented in Appendix \ref{app:est_noise_Y_lemma}.

\begin{lemma}\label{est_noise_Y_lemma}
Let us assume that functions $\phi$ meet Assumption \ref{L_smooth_assump_sp} and for all $j\leq k$ we have $\mathbb{E}\left\|\mathbf{X}^j-\overline{\mathbf{X}}^j\right\|^2\leq \delta'_x,$ $\mathbb{E}\left\|\hat{\mathbf{Y}}^j-\overline{\hat{\mathbf{Y}}}^j\right\|^2 \leq \delta'_y,$ where $\{\mathbf{X}^j\}_{j\leq k}$ and $\{\hat{\mathbf{Y}}^j\}_{j\leq k}$. Then for the bias and noise in the inexact gradient $g_{k,j}=\overline{\widetilde{\nabla}_Y\phi}(\mathbf{X}^k,\mathbf{Y}^k) = \frac{1}{n} \sum\limits_{i=1}^m \widetilde{\nabla_y} \phi_i(x_{k,i}, y_{k,i})$ we  the following estimations:
\begin{equation*}
\mathbb{E}\left\|\mathbb{E}_{X,Y,\xi}g_{k,j}-\nabla_Y \phi(\overline{x}^k, \overline{\hat{y}}^j)\right\|^2\leq \frac{ 3\delta^2 + 3L_{yy}^2 {\delta'_y}+3L_{yx}^2 {\delta'_x}}{n},
\end{equation*}
and
\begin{align*}
\mathbb{E}\left\|\mathbb{E}_{X,Y,\xi}g_{k,j}-g_{k,j}\right\|^2\leq 16\left(\frac{L_{yy,l}^2{\delta'_y}+L_{yx,l}^2\delta'_x+\delta^2}{n}\right)+\frac{18\sigma^2}{n}.
\end{align*}
\end{lemma}

It allows us to use the results of Theorem \ref{theorem:main_theorem_full}. Let parameters $N_y$ and $T_y$ be defined according to this theorem for the current value of outer variable $\mathbf{X}^k$ and chosen accuracies $\varepsilon_y$ and $\delta'_y$. Then after $N_y$ inner iterations we obtain such point $\mathbf{Y} = \hat{\mathbf{Y}}^{N_y}$, s.t.
$$
\max_{y \in \mathbb{R}^{d_y}} \phi(\overline{x}^{k}, y)-\mathbb{E}\phi(\overline{x}^{k}, \overline{\hat{y}}^{N_y})\leq \varepsilon_y+\frac{\Delta_y^2}{2\mu_y n}, \quad \text{and} \quad \mathbb{E}\left\|\hat{\mathbf{Y}}^{N_y}-\overline{\hat{\mathbf{Y}}}^{N_y}\right\|^2 \leq {\delta'_y},
$$
where $\Delta_y^2=19\left(L_{yy,l}^2{\delta'_y}+L_{yx,l}^2\delta'_x+\sigma^2+\delta^2\right).$

This allows us to estimate inaccuracies in the inexact gradient with respect to $\mathbf{X}$ on outer iterations. Note, that it contains new inexactness because of the inexact solution of the inner problem on each iteration. The following lemma contains results about the gradient bias and noise (see the proof in Appendix \ref{app:est_noise_X_lemma}).

\begin{lemma}\label{est_noise_X_lemma}

Let us assume that the function $\phi$ meets Assumption \ref{L_smooth_assump_sp} and for all $j\leq k$ we have $\mathbb{E}\left\|\mathbf{X}^j-\overline{\mathbf{X}}^j\right\|^2 \leq \delta'_x$, $\mathbb{E}\left\|{\mathbf{Y}}^j-\overline{{\mathbf{Y}}}^j \right\|^2 \leq \delta'_y$,  where $\{\mathbf{X}^j\}_j$ and $\{\mathbf{Y}^j\}_j$ are the sequence generated by Algorithm \ref{alg:base_alg_sp}. Also, $\mathbf{Y}^k$ be a point such that $g_x(\overline{y})-g_x^*\leq \varepsilon_y+\frac{\Delta_y^2}{2\mu_y n}$. Then for the bias and noise in the inexact gradient $h_{k} = \overline{\widetilde{\nabla}_{\mathbf{Y}} \phi}(\mathbf{X}^k,\mathbf{Y}^k) = \frac{1}{n}\sum\limits_{i=1}^n \widetilde{\nabla}_y \phi_i(x_{k,i}, y_{k,i})$, we  the following estimations

\begin{equation*}
\mathbb{E}\left\|\mathbb{E}_{\mathbf{X},\mathbf{Y}, \xi} h_{k} - \nabla F(\overline{\mathbf{Y}})\right\|^2\leq \frac{3\delta^2 +3L_{xx}^2 {\delta'_x}+6L_{xy}^2 \left({\frac{2{\varepsilon}_y}{\mu_y}}+\frac{\Delta_y^2 }{\mu_y^2{n}}+{\delta'_y}\right)}{n},
\end{equation*}
and
\begin{align*}
\mathbb{E}\left\|\mathbb{E}_{\mathbf{X}, Y , \xi} h_{k} - h_k\right\|^2 \leq 16 \frac{L_{xy,l}^2{\delta'_y} + L_{xx,l}^2 \delta'_x + \delta^2}{n} + \frac{18\sigma^2}{n}.
\end{align*}
\end{lemma}

Uniting the results of  Lemmas \ref{est_noise_Y_lemma}, \ref{est_noise_X_lemma} and Theorem \ref{theorem:main_theorem_full}, we obtain the following result for the convergence of Algorithm \ref{alg:base_alg_sp}.

\begin{theorem}\label{theorem:main_theorem_sp_full}
Let the functions $\phi(\cdot, y)$ and $-\phi(x, \cdot)$ meet Assumption \ref{PL_assump} with constants $\mu_x>0$ and $\mu_y>0$ and the functions $\phi_i$ meet the Assumption \ref{L_smooth_assump_sp}. Besides statements \eqref{y_minimizer} and \eqref{x_minimizer} hold. The gradient oracle $\left(\widetilde{\nabla}_{\mathbf{X}} \Phi(\mathbf{X},\mathbf{Y}), \widetilde{\nabla}_{\mathbf{Y}} \Phi(\mathbf{X},\mathbf{Y})\right)$ be a $(\delta,\sigma^2)$-biased gradient oracle. Also, Assumption \ref{w_mixing_matrix} holds. Let us introduce the following inexactness values
\begin{align*}
\Delta_y^2&:= 19\left(L_{yy,l}^2{\delta'_y}+L_{yx,l}^2\delta'_x+\sigma^2+\delta^2\right),\\
\Delta_x^2 &:= 22 L_{xy,l}^2{\delta'_y} + 19L_{xx,l}^2\delta'_x+19\delta^2 +18\sigma^2+6L_{xy}^2 \left({\frac{2{\varepsilon}_y}{\mu_y}}+\frac{\Delta_y^2 }{\mu_y^2{n}}\right).
\end{align*}
and define the values $D_{\mathbf{X}}$ and $D_\mathbf{Y}$, such that
\begin{equation*}
D_{\mathbf{X}} = 6\gamma_x^2 \left\| \nabla F(\overline{\mathbf{X}}^*)\right\|^2+2\delta'_x  + 6\left(\gamma_x^2+\frac{1}{\mu_x^2}\right) \Delta^2 +  \frac{12\gamma_x^2 L_x^2}{\mu_x}\left(1-\frac{\mu_x}{L_x}\right) \left(F(\overline{\mathbf{X}}^0)-F^*\right),
\end{equation*}
and $D_{\mathbf{Y}}= \max_k D_{\mathbf{X}^k,\mathbf{Y}}$, where
\begin{align*}
D_{\mathbf{X}^k,\mathbf{Y}} & = 6\gamma^2 \mathbb{E} \left\|\nabla G_{\mathbf{X}^k}(\overline{\mathbf{Y}}^*)\right\|^2 + 2 \delta_y' + 6 \left(\gamma_y^2 + \frac{1}{\mu_y^2} \right) \Delta_y^2 
\\& \;\;\;\; +  \frac{12\gamma^2 L_{yy,g}^2}{\mu_y} \left(1-\frac{\mu_y}{L_{yy,g}}\right) \mathbb{E}\left[G_{\mathbf{X}^k}(\overline{\mathbf{Y}}^0)-G_{\mathbf{X}}^*\right],
\end{align*}
	
Let Algorithm \ref{alg:base_alg_sp}, with step-sizes $\gamma_x=\frac{1}{L_x}, \gamma_y=\frac{1}{L_{yy,g}}$,  makes $ T_x = \tau \left\lceil \frac{1}{2\lambda} \log \left(\frac{D_{\mathbf{X}}}{\delta'_x} \right)\right\rceil$ communication steps in each outer iteration and $T_y = \tau\left\lceil\frac{1}{2\lambda}\log \left(\frac{D_{\mathbf{Y}}}{\delta'_y}\right) \right\rceil$ communication steps in each inner iteration. Then Algorithm \ref{alg:base_alg_sp} requires
$$
N_x= \left\lceil\frac{L_x}{\mu_x}\log \left(\frac{F(\overline{\mathbf{X}}^0)-F^*}{\varepsilon_x}\right)\right\rceil
$$ 
outer iterations, where $L_x=L_{xx,g} + \frac{L_{xy,g}}{\mu_y}$, and
$$
N_y= \left\lceil\frac{L_{yy,g}}{\mu_y}\log \left(\frac{G_{\overline{x}}^*-G_{\overline{x}}(\overline{\mathbf{X}}^0)}{\varepsilon_y}\right) \right\rceil
$$
inner iterations for each outer iteration to yield the pair $(\mathbf{X}^{N_x}, \mathbf{Y}^{N_y})$ such that
\begin{equation}\label{quality_sp_full1}
\mathbb{E}\left[f\left(\overline{x}^{N_x}\right)-f^*\right] \leq \varepsilon_x+\frac{\Delta_x^2}{2\mu_x n},\quad \mathbb{E}\left\|\mathbf{X}^{N_x}-\overline{\mathbf{X}}^{N_x}\right\|\leq \sqrt{\delta'_x},
\end{equation}
and 
\begin{equation}\label{quality_sp_full2}
\mathbb{E}\left[\max_{y \in \mathbb{R}^{d_y}} \phi\left(\overline{x}^{N_x}, y\right)-\phi\left(\overline{x}^{N_x}, \overline{y}^{N_y}\right)\right]\leq \varepsilon_y+\frac{\Delta_y^2}{2\mu_y n}, \quad \mathbb{E}\left\|\mathbf{Y}^{N_y}-\overline{\mathbf{Y}}^{N_y}\right\|\leq \sqrt{\delta'_y}.
\end{equation}
\end{theorem}

\begin{remark}
Note, that the total number of gradient computations and communication steps to approach the qualities \eqref{quality_sp_full1} and \eqref{quality_sp_full2} can be considered comparable with \eqref{quality_sp1} and \eqref{quality_sp2}.
\end{remark}

\begin{remark}
In the case, when $\sigma=0$ (i.e., the noise is almost everywhere is zero), Theorem \ref{theorem:main_theorem_sp_full} gets results worse than Theorem \ref{theorem:main_theorem_sp}. It is related to different rounding and inexact estimations in the proof of Theorem \ref{theorem:main_theorem_sp_full}.
\end{remark}

\section{Numerical experiments}
{\color{black} 
To show the practical performance of the proposed Algorithms \ref{alg:base_alg} and \ref{alg:base_alg_sp}, we performed a series of numerical experiments for the Robust Least Squares problem. All experiments were made using Python 3.4, on a computer with Intel(R) Core(TM) i7-8550U CPU @ 1.80GHz, 1992 Mhz, 4 Core(s), 8 Logical Processor(s), and 8 GB RAM.

Let us consider the following least squares minimization problem
\begin{equation}\label{least_squares}
\min_{x \in \mathbb{R}^{d_x}} \frac{1}{n}\sum\limits_{i=1}^n \frac{1}{2}\left\|A_i x-y_{i,0}\right\|^2,
\end{equation}
for given matrices $A_i\in\mathbb{R}^{d_i\times d_x}$ and vectors $y_{i,0}\in\mathbb{R}^{d_i}$, when $A_i, y_{i,0}$ are placed in $n$ different nodes. In \cite{RobustLS}, it was proposed a robust version of this problem in the following form
\begin{equation*}
\min_{x \in \mathbb{R}^{d_x}} \max_{y\in\mathbb{R}^{d_y}: \|B_i y\|\leq \delta} \frac{1}{n}\sum\limits_{i=1}^n \frac{1}{2}\left\|A_ix - y_{i,0} - B_iy\right\|^2,
\end{equation*}
where $B_i \in\mathbb{R}^{d_i\times d_y}$, for some $\delta >0$. The Robust Least Squares problem with soft constraint has the following form
\begin{equation}\label{robust_least_squares}
\min_{x \in \mathbb{R}^{d_x}} \max_{y\in\mathbb{R}^{d_y}} \left\{\phi(x,y)= \frac{1}{n}\sum\limits_{i=1}^n \phi_i(x, y)\right\},
\end{equation}
where
\begin{equation}\label{functionals_for_exp}
\phi_i(x,y)=\frac{1}{2}\left\|A_ix - y_{i,0} - B_iy\right\|^2 - \frac{\alpha}{2}\|B_i y\|^2,
\end{equation}
for some $\alpha >1$. Note, that the matrices $A_i, B_i$ and the vectors $y_{i,0}$ for each $i \in \{1, \ldots, n\}$, are located in different nodes. So, in the distributed statement, we can not solve this problem explicitly. At the same time, the global solution has a compact form.

Note that the conditions of Theorem \ref{theorem:main_theorem_sp_full} are fulfilled for the considered saddle point problem \eqref{robust_least_squares} with \eqref{functionals_for_exp}. The function $\phi(\cdot, y)$, for each $y \in \mathbb{R}^{d_y}$, satisfies PL-condition with constant $\mu_x$, and  $-\phi(x, \cdot)$, for each $x \in \mathbb{R}^{d_x}$, satisfies PL-condition with constant $\mu_y$, where $\mu_x$ and $\mu_y$ are the smallest non-zero eigenvalues of the matrices $A=\sum\limits_{i=1}^n A_i^\top A_i$ and $B=(\lambda -1)\sum\limits_{i=1}^n B_i^\top B_i$, respectively. Note, that problem \eqref{robust_least_squares} is a convex-concave problem but not strongly-convex-strongly-concave until $A$ and $B$ are not full-rank matrices.

Firstly, for the least squares minimization problem \eqref{least_squares},  we compare the performance of the proposed Algorithm \ref{alg:base_alg} to the recently proposed algorithm DAccGD  \cite{rogozin_accel_do}. In \cite{rogozin_accel_do}, for the least squares problem, it was compared the performance of DAccGD to EXTRA \cite{shi2015extra}, DIGing \cite{nedic2017achieving}, Mudag \cite{ye2023multi}, and APM-C \cite{li2020decentralized}, as a result of the comparison, DAccGD outperformed all the mentioned algorithms.

We run Algorithm \ref{alg:base_alg} and DAccGD, for problem \eqref{least_squares}, with $d_x = 1000$ and  $n = 20$. The matrices $A_i$ and vectors $y_{i,0}$, for each $i \in \{1,\ldots, n\} $,  are randomly generated from the standard normal distribution. The graphs $E^k, \forall k \geq 0$ are randomly generated, and the corresponding mixing matrices $\mathbf{W}^k$ are specified according to \eqref{matrix_W_metropolis}, with Metropolis weights. We take the zero matrices for the initialization. We also choose the fixed step sizes $\gamma_x = 10^{-3}, N_x = 2 \times 10^4, N_y = 100$ and $10$ steps in the consensus at each iteration. For the considered problem \eqref{least_squares}, with the previously mentioned parameters and settings, is ill-conditioned since we have $L \approx 4 \times 10^{6} $ and $\mu \approx 10^{-9}$, thus the conditions number $\kappa = L/\mu \approx 4 \times 10^{15}$.

The results of the conducted experiments are represented in Figure \ref{fig_min_problem}. These results demonstrate the value of the objective function and the norm of its gradient at each generated point $x_k$ by algorithms as a function of iteration $k$. From Fig. \ref{fig_min_problem}, we can see how the proposed Algorithm \ref{alg:base_alg} outperformed DAccGD, and for a not sufficiently big number of iterations (namely $2000$) we can achieve a solution to the problem with high accuracy.

\begin{figure}[htp]
\minipage{0.50\textwidth}
\includegraphics[width=\linewidth]{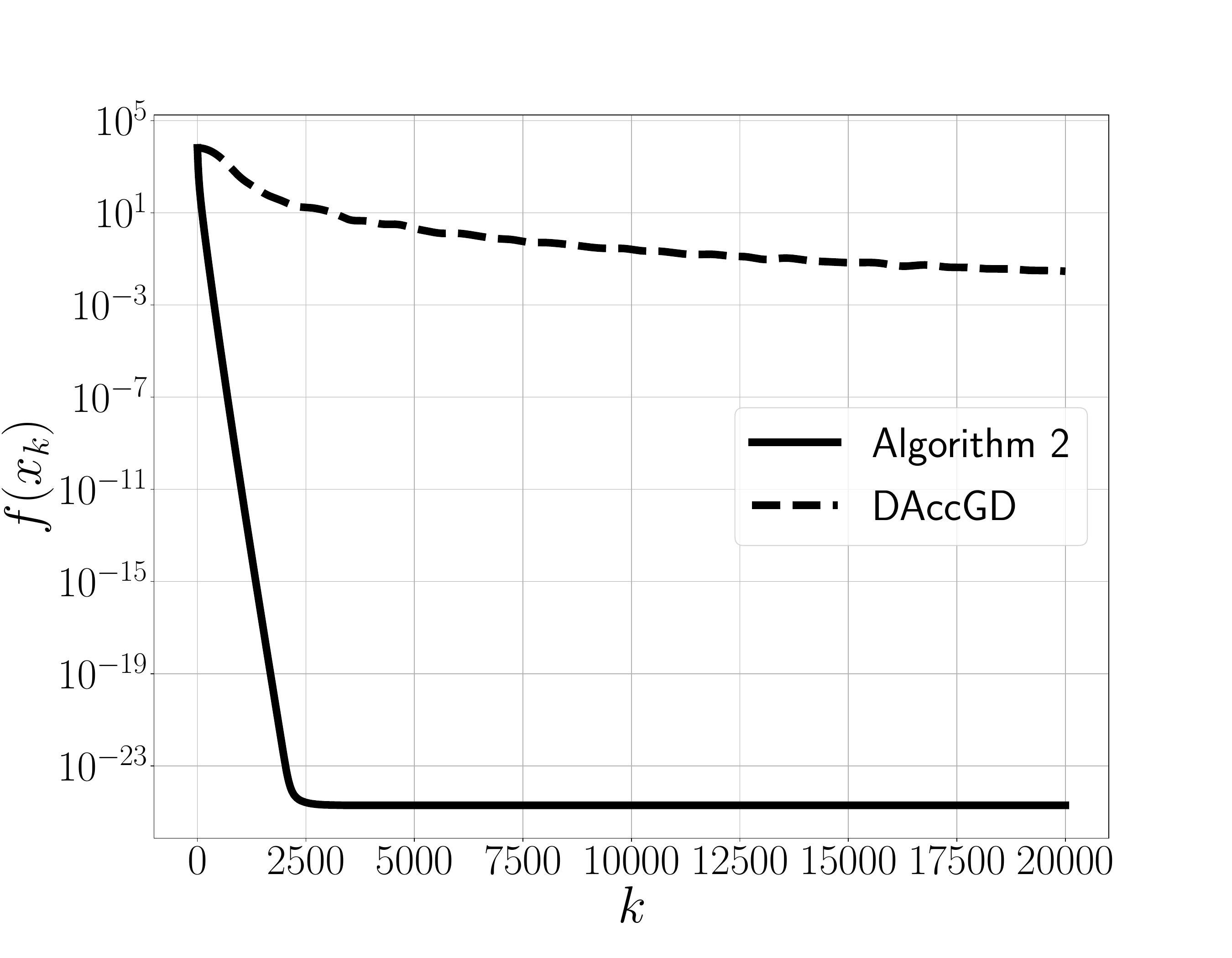}
\endminipage
\minipage{0.50\textwidth}
\includegraphics[width=\linewidth]{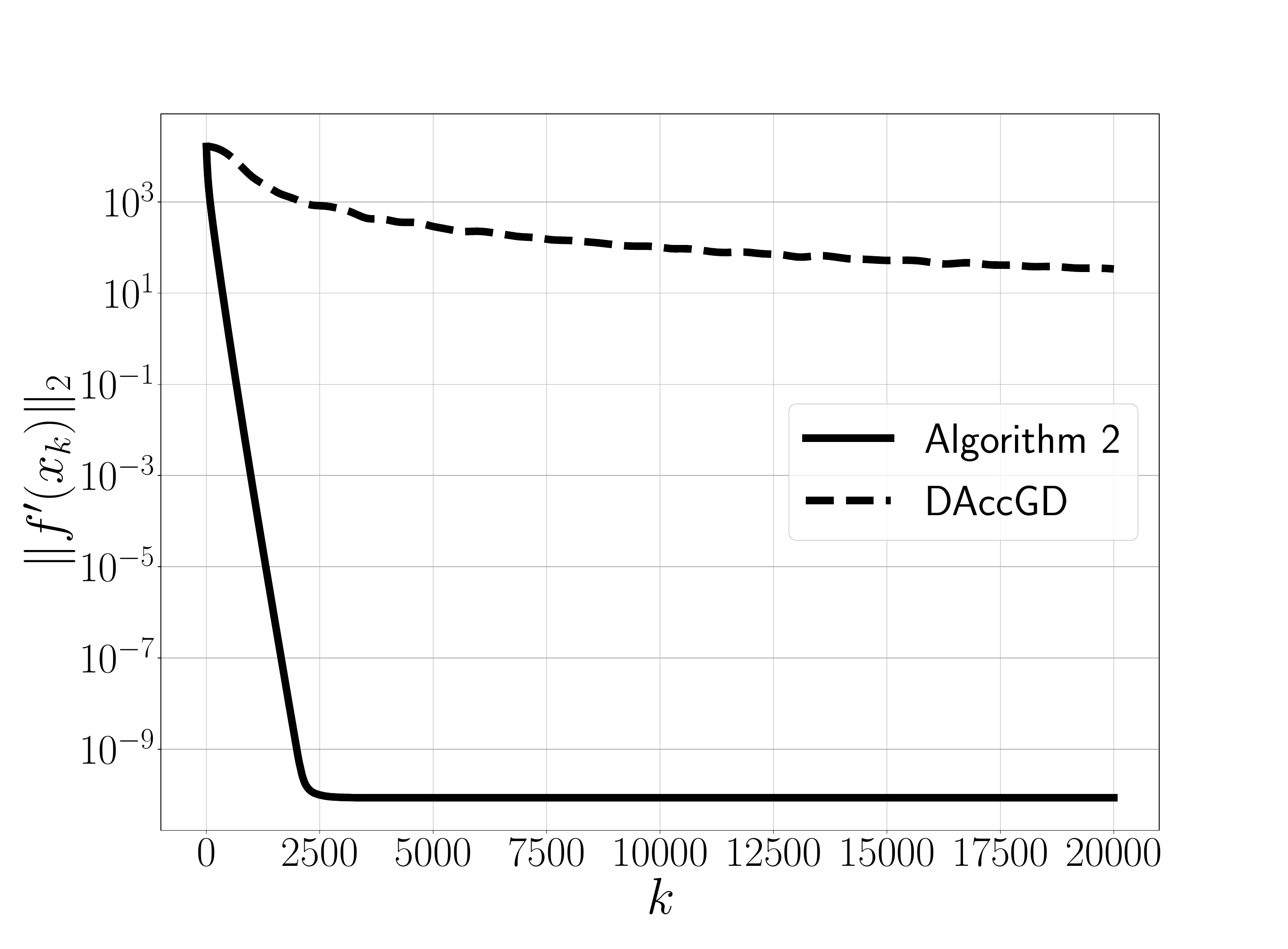}
\endminipage
\caption{Results of  Algorithms \ref{alg:base_alg} and  DAccGD  \cite{rogozin_accel_do}, for  problem  \eqref{least_squares}.}
\label{fig_min_problem}
\end{figure}

Also, for problem  \eqref{robust_least_squares} with \eqref{functionals_for_exp}, we run Algorithm \ref{alg:base_alg_sp} with $d_x = 1000, d_y = 100$ and different values of $n$ (the number of components in \eqref{robust_least_squares}, which indicates the number of nodes in the graph). We take the zero matrices for the initialization of the algorithm. The matrices $A_i, B_i$ and vectors $y_{i,0}$, for each $i \in \{1,\ldots, n\} $ in \eqref{functionals_for_exp} are randomly generated from the standard normal distribution. We also choose the fixed step sizes $\gamma_x = \gamma_y = 10^{-3}, N_x = 10^4, N_y = 10$ (experimentally we see more than $10$ iterations, but there is not a remarkable difference) and $10$ steps in the consensus at each iteration of Algorithm \ref{alg:base_alg_sp}.  The results of the conducted experiments, for problem  \eqref{robust_least_squares} with \eqref{functionals_for_exp},  are represented in Figures \ref{fig_loss_time}, \ref{fig_grad_x_y}. These results demonstrate the value of the objective function $\phi$ and the norm of its gradient with respect to both variables at each point $(x_k, y_k)$ generated by the algorithm as a function of iteration $k$, and the running time of the algorithm in seconds as a function of the number of nodes $n$.

From Fig. \ref{fig_loss_time} and Fig. \ref{fig_grad_x_y}, we can see the efficiency of the proposed Algorithm \ref{alg:base_alg_sp}, which provides a remarkable quality solution with respect to the value of the objective function and the norm of its gradient at each generated point each iteration. Also, we can see the effect of the number of nodes $n$, in the work of the algorithm, where the quality of a solution decreases when we increase (not strongly) the number of nodes.

\begin{figure}[htp]
\minipage{0.50\textwidth}
\includegraphics[width=\linewidth]{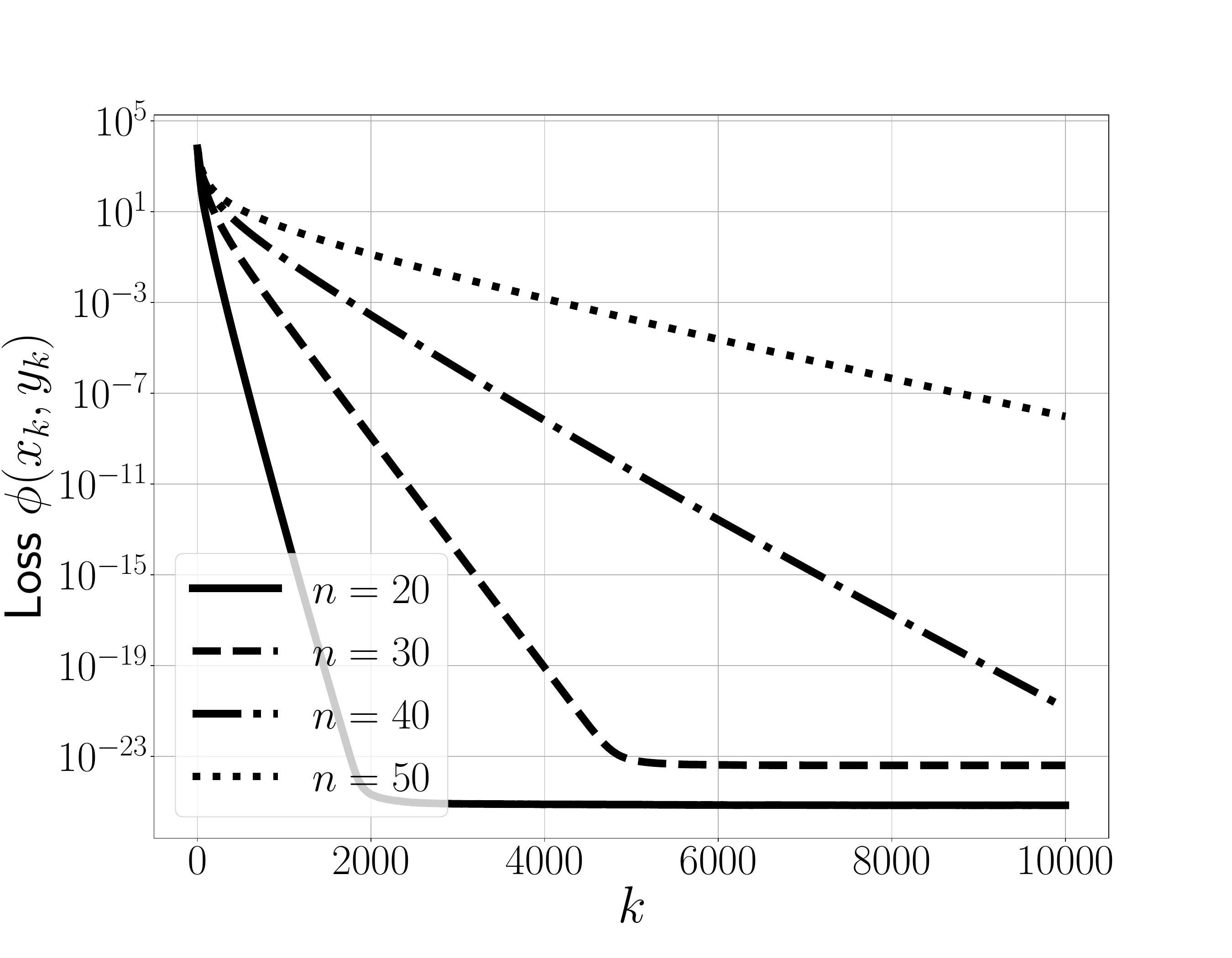}
\endminipage
\minipage{0.50\textwidth}
\includegraphics[width=\linewidth]{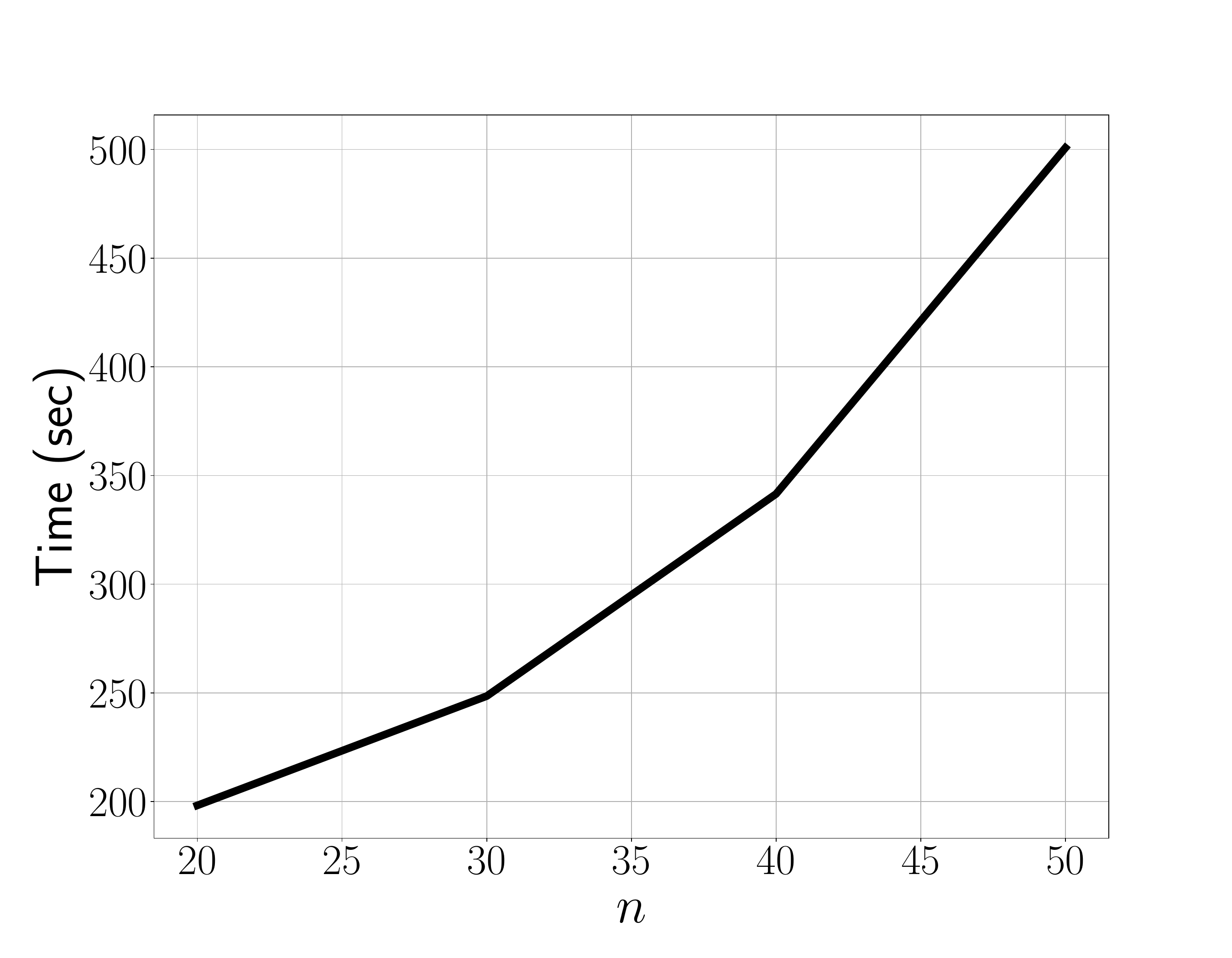}
\endminipage
\caption{Results of  Algorithm \ref{alg:base_alg_sp} for  problem  \eqref{robust_least_squares} with \eqref{functionals_for_exp}, $\alpha = 2$, and different values of $n$.}
\label{fig_loss_time}
\end{figure} 

\begin{figure}[htp]
\minipage{0.50\textwidth}
\includegraphics[width=\linewidth]{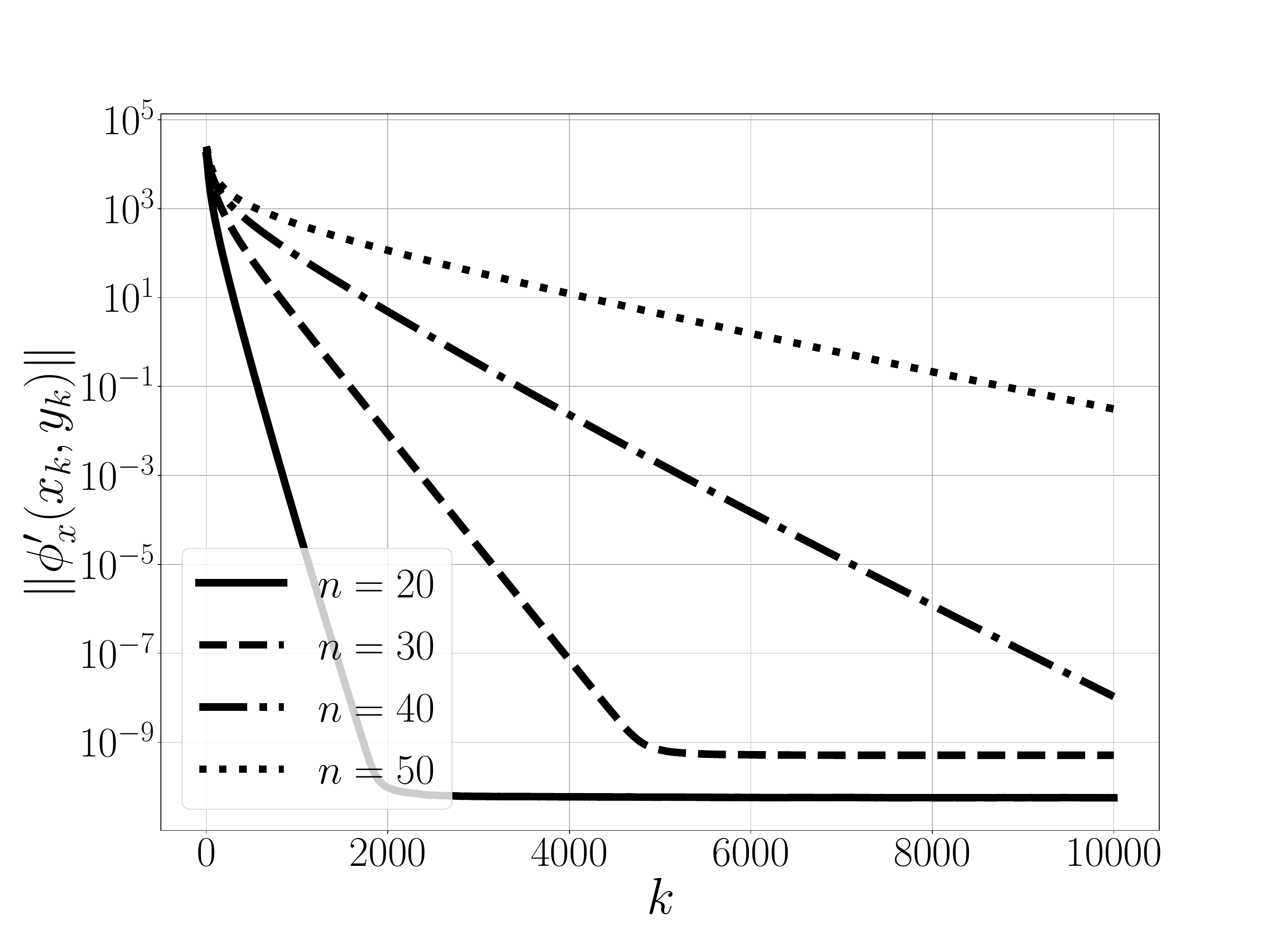}
\endminipage
\minipage{0.50\textwidth}
\includegraphics[width=\linewidth]{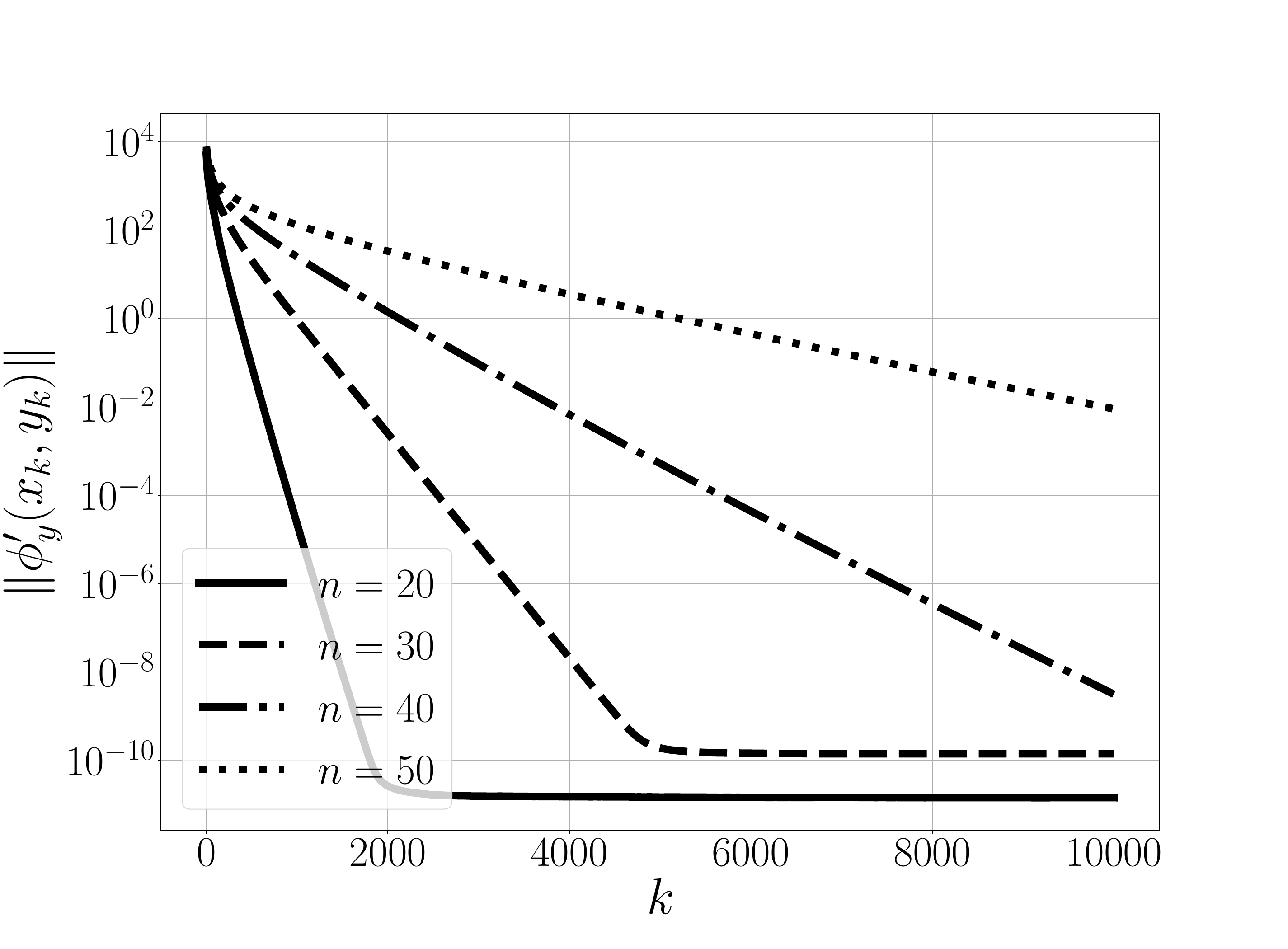}
\endminipage
\caption{Results of  Algorithm \ref{alg:base_alg_sp} for problem  \eqref{robust_least_squares} with \eqref{functionals_for_exp}, $\alpha = 2$, and  different values of $n$.}
\label{fig_grad_x_y}
\end{figure} 
}




\section{Conclusion}

In this paper, we studied the decentralized minimization problem with objective functions that satisfy Polyak-{\L}ojasiewicz condition (PL-condition), and the decentralized saddle-point problem with a structure of the sum, with objectives satisfying the two-sided PL-condition. For solving the minimization problem under consideration, we proposed a gradient descent type method with a consensus projection procedure and the inexact gradient of the objective function. To solve the considered saddle point problem,  we proposed a generalization of the Multi-step Gradient Descent Ascent method with a consensus procedure, and inexact gradients of the objective function concerning both variables. For the studied classes of the problems, we estimated the sufficient communication steps and iterations number to approach the required quality concerning the function and the distance between the agent's vectors and their average some results of the conducted numerical experiments are presented, which demonstrate the effectiveness of the proposed algorithm for the least squares and robust least squares problems.

\bibliography{biblio_new}

\appendix

\section{Proof of Lemma \ref{est_noise_lemma}.} \label{app:noise_est}
\begin{proof}
Note, for the bias in $\overline{\widetilde{\nabla}} f(\overline{\mathbf{X}}^k)$ we proved in the main part the following estimation
\begin{equation}
\mathbb{E}\left\|\mathbb{E}_{\mathbf{X}^k,\xi} \overline{\widetilde{\nabla} f}(\mathbf{X}^k)-{{\nabla} f}(\overline{x}^k)\right\|^2\leq \frac{2\delta^2 + L_l^2 \delta'}{{n}}.
\end{equation}
Also, we assume that the condition $\mathbb{E}_{\mathbf{X}^j} \left\|\mathbf{X}^j-\overline{\mathbf{X}}^j\right\|^2 \leq \delta'$ holds for all $j\leq k$.

Let us estimate noise in $\overline{\widetilde{\nabla}} f(\overline{\mathbf{X}}^k)$. For this, we estimate the second moment of the random component in $\widetilde{\nabla} F(\mathbf{X}^k)$ for given $\overline{\mathbf{X}}^k$. It is given by the expression
$$
\mathbb{E}_{\mathbf{X}^k,\xi}\left\|\mathbb{E}_{\mathbf{X}^k,\xi}\widetilde{\nabla} F(\mathbf{X}^k)-\widetilde{\nabla} F(\mathbf{X}^k)\right\|.
$$

Let us estimate the inner part in the following way:
$$
\begin{aligned}
\left\|\mathbb{E}_{\mathbf{X}^k,\xi}\widetilde{\nabla} F(\mathbf{X}^k)-\widetilde{\nabla} F(\mathbf{X}^k)\right\| & \leq \left\|\mathbb{E}_{\mathbf{X}^k,\xi}\widetilde{\nabla} F(\mathbf{X}^k)-\mathbb{E}_{\mathbf{X}^k}\widetilde{\nabla} F(\mathbf{X}^k)\right\| \\& \quad +\left\|\mathbb{E}_{\mathbf{X}^k}\widetilde{\nabla} F(\mathbf{X}^k)-\widetilde{\nabla} F(\mathbf{X}^k)\right\|
\\& \leq \mathbb{E}_{\mathbf{X}^k}\left\|\mathbb{E}_{\xi}\widetilde{\nabla} F(\mathbf{X}^k)-\widetilde{\nabla} F(\mathbf{X}^k)\right\| 
\\& \quad +\left\| \mathbb{E}_{\mathbf{X}^k} \widetilde{\nabla} F(\mathbf{X}^k) -\widetilde{\nabla} F(\mathbf{X}^k) \right\|,
\end{aligned}
$$
\rev{where we used triangle inequality and Jensen's inequality. }

\rev{Using these statements again, the second term in the sum above can be represented in the following way:}
\begin{align*}
\left\|\mathbb{E}_{\mathbf{X}^k}\widetilde{\nabla} F(\mathbf{X}^k)-\widetilde{\nabla} F(\mathbf{X}^k)\right\|\leq& \mathbb{E}_{\mathbf{X}^k}\left\|\widetilde{\nabla} F(\mathbf{X}^k)-{\nabla} F(\mathbf{X}^k)\right\| 
\\&+
\left\|\mathbb{E}_{\mathbf{X}^k}{\nabla} F(\mathbf{X}^k)-{\nabla} F(\mathbf{X}^k)\right\| + \left\|\widetilde{\nabla} F(\mathbf{X}^k)-{\nabla} F(\mathbf{X}^k)\right\|.
\end{align*}
From this, we have the following estimation for its square expectation:
$$
\begin{aligned}
\mathbb{E} \left\|\mathbb{E}_{\mathbf{X}^k}\widetilde{\nabla} F(\mathbf{X}^k)-\widetilde{\nabla} F(\mathbf{X}^k)\right\|^2 \leq &  2\mathbb{E}\left\|\mathbb{E}_{\mathbf{X}^k}{\nabla} F(\mathbf{X}^k)-{\nabla} F(\mathbf{X}^k)\right\|^2 
\\& + 8\mathbb{E}\left\|\widetilde{\nabla} F(\mathbf{X}^k)-{\nabla} F(\mathbf{X}^k)\right\|^2.
\end{aligned}
$$

The  term $\left\|\mathbb{E}_{\mathbf{X}^k}{\nabla} F(\mathbf{X}^k)-{\nabla} F(\mathbf{X}^k)\right\|$  can be estimated above by the following sum
$$
\mathbb{E}_{\mathbf{X}^k}\left\|{\nabla} F(\mathbf{X}^k)-{\nabla} F(\overline{\mathbf{X}}^k)\right\|+\left\|{\nabla} F(\overline{\mathbf{X}}^k)-{\nabla} F(\mathbf{X}^k)\right\|.
$$
From this we can estimate the full mathematical expectation as $\mathbb{E}\left\|\mathbb{E}_{\mathbf{X}^k}{\nabla} F(\mathbf{X}^k)-{\nabla} F(\mathbf{X}^k)\right\|^2 \leq 4L_l^2 \delta'$ until $\mathbb{E}\left\|\mathbf{X}^k-\overline{\mathbf{X}}^k\right\|^2\leq \delta'$. Let us introduce the noise value $n(\mathbf{X}^k)=\widetilde{\nabla} F(\mathbf{X}^k)-\mathbb{E}_{\xi}\widetilde{\nabla} F(\mathbf{X}^k)$ and the bias $b(\mathbf{X}^k)=\mathbb{E}_{\xi}\widetilde{\nabla} F(\mathbf{X}^k)-{\nabla} F(\mathbf{X}^k).$ Note, that $\mathbb{E}\|n(\mathbf{X}^k)\|^2\leq \sigma^2$ and $\mathbb{E}\left\|b(\mathbf{X}^k)\right\|^2 \leq \delta^2$. So, the expectation of another term can be estimated in the following form: 
\[
\mathbb{E} \left\|\widetilde{\nabla} F(\mathbf{X}^k)-{\nabla} F(\mathbf{X}^k)\right\|^2=\mathbb{E}\left\|b(\mathbf{X}^k)\right\|^2 + \mathbb{E}\left\|n(\mathbf{X}^k)\right\|^2\leq \sigma^2+\delta^2.
\]
Finally, we have that
\[
\mathbb{E}\left\|\mathbb{E}_{\mathbf{X}^k}\widetilde{\nabla} F(\mathbf{X}^k)-\widetilde{\nabla} F(\mathbf{X}^k)\right\|^2 \leq 8L_l\sqrt{\delta'}+8\sigma+8\delta.
\] 

Uniting the inequalities above, we have the following estimation:
\begin{align*}
\mathbb{E}\left\|\mathbb{E}_{\mathbf{X}^k,\xi}\widetilde{\nabla} F(\mathbf{X}^k)-\widetilde{\nabla} F(\mathbf{X}^k)\right\|^2 \leq &  2 \mathbb{E}_{\mathbf{X}^k,\xi}\left(\mathbb{E}_{\mathbf{X}^k}\left\|\mathbb{E}_{\xi}\widetilde{\nabla} F(\mathbf{X}^k)-\widetilde{\nabla} F(\mathbf{X}^k)\right\|\right)^2
\\&
+ 16L_l^2{\delta'}+16\sigma^2+16\delta^2\\
\leq&16L_l^2{\delta'}+18\sigma^2+16\delta^2.
\end{align*}

Finally, for the noise component in $\overline{\widetilde{\nabla}} f(\overline{X}^k)$, we have the following estimation:
\begin{align*}
\mathbb{E}\left\|\mathbb{E}_{X,\xi}\overline{\widetilde{\nabla} f}(X^k)-\overline{\widetilde{\nabla} f}(X^k)\right\|^2= & \frac{1}{n}\mathbb{E}\left\|\mathbb{E}_{\mathbf{X}^k,\xi}\overline{\widetilde{\nabla} F}(\mathbf{X}^k)-\overline{\widetilde{\nabla} F}(\mathbf{X}^k)\right\|^2
\\   \leq&
\frac{1}{n}\mathbb{E}\left\|\mathbb{E}_{\mathbf{X}^k,\xi}{\widetilde{\nabla} F}(\mathbf{X}^k)-{\widetilde{\nabla} F}(\mathbf{X}^k)\right\|^2
\\\leq& 
\frac{16L_l^2{\delta'}+18\sigma^2+16\delta^2}{n}.
\end{align*}
So, we obtained the estimation for the noise component in $\overline{\widetilde{\nabla} f}(X^k)$.
\qed

\section{Proof of Lemma \ref{lemma_conv}.}\label{app:proof_lemma_conv}
\begin{proof}
According to proof of Lemma 2 from \cite{stich_sgd_biased} we can prove the following inequality:
\[
\mathbb{E}\left[f(x^{t+1})\Big| x^t\right] \leq f(x^t) - \frac{\gamma}{2} \left\|\nabla f(x^t)\right\|^2 + \frac{\gamma}{2} \left\|b(x^k)\right\|^2 +\frac{\gamma^2 L}{2}\mathbb{E}_\xi \left\|n(x^k)\right\|^2.
\]
Using PL-condition for the function $f$ and taking full mathematical expectation, we have the following inequality
\[
\mathbb{E}\left[f(x^{t+1})-f^*\right]\leq (1-\gamma\mu)\mathbb{E}\left[f(x^t)-f^*\right] + \frac{\gamma}{2}\mathbb{E} \left \|b(x^k)\right\|^2 +\frac{\gamma^2 L}{2}\mathbb{E}\|n(x^k)\|^2.
\]
Using lemma's conditions we obtain the following inequality:
\[
\mathbb{E}\left[f(x^{t+1})-f^*\right]\leq (1-\gamma\mu)\mathbb{E}\left[f(x^t)-f^*\right] + \frac{\gamma \delta^2 }{2} +\frac{\gamma^2 L\sigma^2}{2},
\]
or
\[
\mathbb{E}\left[f(x^k)-f^*\right]\leq (1-\gamma\mu)^{k}\left(f(x^0)-f^*\right) + \frac{ \delta^2+\gamma L\sigma^2 }{2\mu}.
\]
\end{proof}

\section{Proof of Theorem \ref{theorem:main_theorem_full}.} \label{app:proof_main_theorem_full}
\begin{proof}
Using Assumptions \ref{PL_assump}, Lemma \ref{lemma_conv}, taking step size $\gamma=\frac{1}{L_g}$, we can estimate convergence for function $f$:
\begin{equation}
\mathbb{E}\left[f(\overline{x}^{k})-f^*\right]\leq \left(1-\frac{\mu}{L_g}\right)^{k} (f(\overline{x}^0)-f^*) + \frac{\Delta^2}{2\mu n},
\end{equation}
until $\mathbb{E} \left\|\mathbf{X}^j-\overline{\mathbf{X}}^j\right\| \leq\delta^{\prime}$ for $j< k$ where $\Delta^2=18\left(L_l^2{\delta^{\prime}}+\sigma^2+\delta^2\right)$.

But really on iterations of Algorithm \ref{alg:base_alg} we have access only to $\mathbf{X}^k$. Let us find the sufficient communication step such that the following expression is true:
\begin{equation}\label{rand_cond_delta}
\mathbb{E} \left\|\mathbf{X}^{k}-\overline{\mathbf{X}}^{k}\right\|^2 \leq {\delta^{\prime}} \Longrightarrow \mathbb{E}\left\|\mathbf{X}^{k+1}-\overline{\mathbf{X}}^{k+1}\right\|^2\leq {\delta^{\prime}}.
\end{equation}

By contraction property, we have
\begin{equation}\label{ineq_000}
\left\|\mathbf{X}^{k+1}-\overline{\mathbf{X}}^{k+1}\right\|\leq (1-\lambda)^{\left\lfloor\frac{T_k}{\tau}\right\rfloor} \left\|\mathbf{Z}^{k+1}-\overline{\mathbf{X}}^{k+1}\right\|.
\end{equation}
Here we used, that $\overline{\mathbf{Z}}^{k+1}=\overline{\mathbf{X}}^{k+1}$. Let us estimate the right part of \eqref{ineq_000}.
\begin{align*}
\mathbb{E}\left\|\overline{\mathbf{X}}^{k+1}-\mathbf{Z}^{k+1}\right\|^2 
&
\leq 2\mathbb{E} \left\|\overline{\mathbf{X}}^k-\mathbf{X}^k \right\|^2 + 2\gamma^2 \mathbb{E} \left\| \widetilde{\nabla}F(\mathbf{X}^k) \right\|^2
\\&
\leq 2\delta^{\prime} + 6\gamma^2 \mathbb{E} \left\| \widetilde{\nabla} F(\mathbf{X}^k) - \nabla F(\overline{\mathbf{X}}^k) \right\|^2
\\& \quad 
+6\gamma^2 \mathbb{E} \left \|\nabla F(\overline{\mathbf{X}}^k)- \nabla F(\overline{\mathbf{X}}^*)\right\|^2 + 6\gamma^2 \left\|\nabla F(\overline{\mathbf{X}}^*)\right\|^2 
\\&
\leq 2\delta'+6\gamma^2 \Delta^2 +6\gamma^2 L_g^2  \mathbb{E} \left\|\overline{\mathbf{X}}^k-\overline{\mathbf{X}}^*\right\|^2+6\gamma^2 \left\|\nabla F(\overline{\mathbf{X}}^*)\right\|^2.
\end{align*}

From quadratic growth condition \eqref{QC_cond} we have 
$$
\mathbb{E} \left\|\overline{\mathbf{X}}^k - \overline{\mathbf{X}}^*\right\|^2 
= n\mathbb{E} \left\|\overline{x}^k - \overline{x}^*\right\|^2 \leq \frac{2n}{\mu} \mathbb{E} \left[f(x^k)-f^*\right].
$$ 
Further, using the convergence rate \eqref{convergence_rate_lemma} we can estimate the third value in the sum above:
$$
\mathbb{E} \left\|\overline{\mathbf{X}}^k-\overline{\mathbf{X}}^*\right\|^2 \leq \frac{2n}{\mu}\left(1-\frac{\mu}{L_g}\right)^{k+1} \left(F(\overline{\mathbf{X}}^0)-F^*\right) + \frac{\Delta^2}{\mu^2}.
$$
Note, that $\mathbb{E}\xi\leq \sqrt{\mathbb{E}\xi^2}$ for any random value $\xi$. Finally, we have that:
$$
\mathbb{E} \left\|\overline{\mathbf{X}}^{k+1}-\mathbf{Z}^{k+1}\right\|^2 \leq D,
$$
for $D$ such that:
$$
D=6\gamma^2 \|\nabla F(\overline{\mathbf{X}}^*)\|^2+2\delta'+6\left(\gamma^2+\frac{1}{\mu^2}\right) \Delta^2 +  \frac{12\gamma^2 L_g^2}{\mu}\left(1-\frac{\mu}{L_g}\right) \left(F(\overline{\mathbf{X}}^0)-F^*\right).
$$

So, for $T_k=\tau\left\lceil\frac{1}{2\lambda}\log \left(\frac{D}{\delta'}\right)\right\rceil$ the condition \eqref{rand_cond_delta} is met.
\end{proof}

\section{Proof of Lemma \ref{est_noise_Y_lemma}.} \label{app:est_noise_Y_lemma}
{\it Proof} For the given $\overline{\mathbf{X}}$ and $\overline{\mathbf{Y}}$, we can estimate the bias in inexact gradient with respect to $Y$ in the following form:
\begin{align*}
\mathbb{E}\left\|\mathbb{E}_{\mathbf{X}^k,\mathbf{Y},\xi}\widetilde{\nabla}_{\mathbf{Y}}\Phi(\mathbf{X}, \mathbf{Y})-\nabla_Y \Phi(\overline{\mathbf{X}}, \overline{\mathbf{Y}})\right\|^2 \leq &
3\mathbb{E} \left\|\mathbb{E}_{\mathbf{X}^k,\mathbf{Y},\xi}\widetilde{\nabla}_{\mathbf{Y}}\Phi(\mathbf{X},\mathbf{Y})-\mathbb{E}_{\mathbf{X}^k,\mathbf{Y}}\nabla_{\mathbf{Y}} \Phi(\mathbf{X}, \mathbf{Y}) \right\|^2 
\\&
+3\mathbb{E}\left\| \mathbb{E}_{\mathbf{X}^k, \mathbf{Y}}{\nabla}_Y\Phi(\mathbf{X}, \mathbf{Y})-\mathbb{E}_{\mathbf{X}^k}\nabla_{\mathbf{Y}} \Phi(\mathbf{X}, \overline{\mathbf{Y}}) \right\|^2 
\\&
+3\mathbb{E} \left\| \mathbb{E}_{\mathbf{X}^k}{\nabla}_{\mathbf{Y}}\Phi(\mathbf{X}, \overline{\mathbf{Y}})-\nabla_{\mathbf{Y}} \Phi(\overline{\mathbf{X}}, \overline{\mathbf{Y}}) \right\|^2
\\ \leq & 
3\mathbb{E} \left\|\mathbb{E}_{\xi}\widetilde{\nabla}_Y\Phi(\mathbf{X},\mathbf{Y})-\nabla_Y \Phi(\mathbf{X}, \mathbf{Y})\right\|^2
\\&
+3\mathbb{E}\|{\nabla}_{\mathbf{Y}}\Phi(\mathbf{X}, \mathbf{Y})-\nabla_Y \Phi(\mathbf{X}, \overline{\mathbf{Y}})\|^2
\\&
+3\mathbb{E} \left\|{\nabla}_Y\Phi(\mathbf{Y}, \overline{\mathbf{Y}})-\nabla_{\mathbf{Y}} \Phi(\overline{\mathbf{X}}, \overline{\mathbf{Y}})\right\|^2
\\ \leq&
3\delta^2 + 3L_{yy}^2 {\delta^{\prime}_y}+3L_{yx}^2 {\delta^{\prime}_x}.
\end{align*}
In the inequality above we used fact that $\|a+b+c\|^2\leq 3\|a\|^2+3\|b\|^2+3\|c\|^2$ for each vectors $a, b$ and $c$. Note, for the bias in $\overline{\widetilde{\nabla}} f(\overline{X}^k)$ we proved in the main part the following estimation
\begin{equation}
\mathbb{E} \left\|\mathbb{E}_{X,\xi} \overline{\widetilde{\nabla} f}(X^k)-{{\nabla} f}(\overline{x}^k) \right\|^2 \leq \frac{2\delta^2 + L_l^2 \delta^{\prime}}{{n}}.
\end{equation}
Also, we assume that the condition $\mathbb{E}_{\mathbf{X}^k} \left\|\mathbf{X}^j-\overline{\mathbf{X}}^j\right\|^2 \leq \delta^{\prime}$ holds for all $j\leq k$.

Let us estimate noise in $\overline{\widetilde{\nabla}} f(\overline{X}^k)$. For this, we estimate the second moment of the random component in $\widetilde{\nabla} F(\mathbf{X}^k)$ for given $\overline{\mathbf{X}}^k$. It is given by the expression $\mathbb{E}_{\mathbf{X}^k,\xi}\|\mathbb{E}_{\mathbf{X}^k,\xi}\widetilde{\nabla} F(\mathbf{X}^k)-\widetilde{\nabla} F(\mathbf{X}^k)\|$. Let us estimate the inner part in the following way:
\begin{align*}
\left\|\mathbb{E}_{\mathbf{X},\mathbf{Y},\xi}\widetilde{\nabla}_{\mathbf{Y}} \Phi(\mathbf{X},\mathbf{Y})-\widetilde{\nabla}_Y \Phi(\mathbf{X}, \mathbf{Y})\right\| 
\leq& 
\left\|\mathbb{E}_{\mathbf{X},\mathbf{Y},\xi}\widetilde{\nabla}_{\mathbf{Y}} \Phi(\mathbf{X},\mathbf{Y})-\mathbb{E}_{\mathbf{X},\mathbf{Y}}\widetilde{\nabla}_{\mathbf{Y}} \Phi(\mathbf{X}, \mathbf{Y})\right\|
\\&+
\left\|\mathbb{E}_{\mathbf{X},\mathbf{Y}}\widetilde{\nabla}_Y \Phi(\mathbf{X}, \mathbf{Y})-\widetilde{\nabla}_Y \Phi(\mathbf{X}, \mathbf{Y})\right\|
\\\leq&
\mathbb{E}_{\mathbf{X},\mathbf{Y}}
\left\|\mathbb{E}_{\xi}\widetilde{\nabla}_Y \Phi(X,Y)-\widetilde{\nabla}_{\mathbf{Y}} \Phi(\mathbf{X}, \mathbf{Y})\right\|
\\&+
\left\|\mathbb{E}_{\mathbf{X},\mathbf{Y}}\widetilde{\nabla}_{\mathbf{Y}} \Phi(\mathbf{X}, \mathbf{Y}) - \widetilde{\nabla}_{\mathbf{Y}}\Phi(\mathbf{X}, \mathbf{Y})\right\|.
\end{align*}
The last term can be estimated in a similar way as in Appendix \ref{app:noise_est}:
\begin{align*}
\mathbb{E} \left\|\mathbb{E}_{\mathbf{X},\mathbf{Y}}\widetilde{\nabla}_{\mathbf{Y}} \Phi(\mathbf{X}, \mathbf{Y}) - \widetilde{\nabla}_{\mathbf{Y}} \Phi(\mathbf{X}, \mathbf{Y})\right\|^2 \leq 
&
2\mathbb{E} \left\|\mathbb{E}_{\mathbf{X},\mathbf{Y}}{\nabla}_{\mathbf{Y}} \Phi(\mathbf{X}, \mathbf{Y})-{\nabla}_{\mathbf{Y}} \Phi(\mathbf{X}, \mathbf{Y})\right\|^2 
\\&
+ 8 \sigma^2 + 8\delta^2
\\ \leq & 
8L_{yy,l}^2{\delta^{\prime}_y} + 8 L_{yx,l}^2 \delta^{\prime}_x + 8\sigma^2 + 8 \delta^2,
\end{align*}

Finally, we can estimate required the second moment of noise component:
\begin{align*}
\mathbb{E} \left\|\mathbb{E}_{\mathbf{X},\mathbf{Y},\xi}\widetilde{\nabla}_{\mathbf{Y}} \Phi(\mathbf{X},\mathbf{Y})-\widetilde{\nabla}_{\mathbf{Y}} \Phi(\mathbf{X}, \mathbf{Y})\right\|^2\leq
&
2\mathbb{E} \left\|\mathbb{E}_{\xi}\widetilde{\nabla}_Y \Phi(\mathbf{X},\mathbf{Y})-\widetilde{\nabla}_{\mathbf{Y}} \Phi(\mathbf{X}, \mathbf{Y})\right\|^2
\\& 
+ 16 L_{yy,l}^2 {\delta^{\prime}_y} + 16 L_{yx,l}^2 \delta^{\prime}_x + 16 \sigma^2 + 16 \delta^2
\\ \leq & 
16L_{yy,l}^2{\delta^{\prime}_y} + 16L_{yx,l}^2 \delta^{\prime}_x + 18\sigma^2 + 16 \delta^2. 
\end{align*}
Thus we can conclude, that 
\begin{equation}\label{est_noise_y_app}
\mathbb{E} \left\|\mathbb{E}\overline{\widetilde{\nabla}_{\mathbf{Y}} \phi}(\mathbf{X},\mathbf{Y})-\overline{\widetilde{\nabla}_{\mathbf{Y}} \phi}(\mathbf{X},\mathbf{Y})\right\| \leq \frac{16L_{yy,l}^2{\delta^{\prime}_y} + 16 L_{yx,l}^2\delta^{\prime}_x + 18\sigma^2+16\delta^2}{n}.
\end{equation}
\end{proof}

\section{ Proof of Lemma \ref{est_noise_X_lemma}} \label{app:est_noise_X_lemma}
\begin{proof}
In a similar way, like for the gradient with respect to the inner variable $\mathbf{Y}$, we can estimate the mathematical expectation of bias with respect to the outer variable $\mathbf{X}$:
$$
\begin{aligned}
\mathbb{E} \left\|\mathbb{E}_{\mathbf{X},\mathbf{Y},\xi}\widetilde{\nabla}_{\mathbf{X}} \Phi(\mathbf{X},\mathbf{Y})-\nabla F\left(\overline{\mathbf{X}}\right)\right\|^2 & \leq
3 \mathbb{E} \left\|\mathbb{E}_{\mathbf{X},\mathbf{Y},\xi}\widetilde{\nabla}_{\mathbf{X}} \Phi(\mathbf{X},\mathbf{Y}) - \mathbb{E}_{\mathbf{X},\mathbf{Y}} \nabla_{\mathbf{X}} \Phi(\mathbf{X}, \mathbf{Y})\right\|^2
\\ & 
+3 \mathbb{E} \left\|\mathbb{E}_{\mathbf{X}, Y}{\nabla}_{\mathbf{X}} \Phi(\mathbf{X}, \mathbf{Y})-\mathbb{E}_{\mathbf{X}}\nabla_{\mathbf{X}} \Phi(\mathbf{X}, \overline{\mathbf{Y}}^*(\mathbf{X}))\right\|^2
\\ & 
+3\mathbb{E} \left\|\mathbb{E}_{\mathbf{X}}{\nabla}_{\mathbf{X}} \Phi(\mathbf{X}, \overline{\mathbf{Y}}^*(\mathbf{X}))-\nabla_{\mathbf{X}} \Phi(\overline{\mathbf{X}}, \overline{\mathbf{Y}}^*(\mathbf{X}))\right\|^2
\\  & \leq
3 \mathbb{E} \left\|\mathbb{E}_{\xi}\widetilde{\nabla}_{\mathbf{X}} \Phi(\mathbf{X},\mathbf{Y}) - \nabla_{\mathbf{X}} \Phi(\mathbf{X}, \mathbf{Y})\right\|^2
\\ &
+ 3\mathbb{E} \left\|{\nabla}_{\mathbf{X}} \Phi(\mathbf{X}, \mathbf{Y})-\nabla_{\mathbf{X}} \Phi (\mathbf{X}, \overline{\mathbf{Y}}^* (\mathbf{X}))\right\|^2
\\ &
+3\mathbb{E} \left\|{\nabla}_{\mathbf{X}} \Phi(\mathbf{X}, \overline{\mathbf{Y}}^*(\mathbf{X})-\nabla_{\mathbf{X}} \Phi(\overline{\mathbf{X}}, \overline{\mathbf{Y}}^*(\mathbf{X}))\right\|^2
\\  &  \leq
3 \delta^2 +3L_{xx}^2 {\delta'_x}+6L_{xy}^2 \left({\frac{2{\varepsilon}_y}{\mu_y}}+\frac{\Delta_y^2 }{\mu_y^2{n}}+{\delta'_y}\right).    
\end{aligned}
$$

Note, that the noise component in the inexact gradient $\overline{\widetilde{\nabla}_{\mathbf{X}} \phi}(\mathbf{X},\mathbf{Y})$  does not depend on inexactness on a solution by the inner problem. So, repeating the steps above for inexact gradient with respect to the variable $\mathbf{X}$, we can obtain the following estimation
\begin{equation}\label{est_noise_x^app}
\mathbb{E}\left\|\mathbb{E}\overline{\widetilde{\nabla}_{\mathbf{X}} \phi}(\mathbf{X},\mathbf{Y})-\overline{\widetilde{\nabla}_{\mathbf{X}} \phi}(\mathbf{X},\mathbf{Y})\right\| \leq \frac{16L_{xy,l}^2{\delta'_y}+16L_{xx,l}^2\delta'_x+18\sigma^2+16\delta^2}{n}.
\end{equation}

Note, that the variable $Y$ in \eqref{est_noise_x^app} is a result of the inner loop in Algorithm \ref{alg:base_alg_sp}. The estimations \eqref{est_noise_y_app} and \eqref{est_noise_x^app} give the bounds on the second norm of noise.
\end{proof}

\end{document}